\documentclass[12pt]{article}

\usepackage{geometry}
\usepackage{amsmath}
\usepackage{amsthm}
\usepackage{amsfonts,dsfont}
\usepackage{amssymb}
\usepackage{stmaryrd}
\usepackage[all,cmtip]{xy}
\usepackage{fancyhdr}
\usepackage{verbatim}
\usepackage{color}
\usepackage{graphicx}
\usepackage[lofdepth,lotdepth]{subfig}
\usepackage{psfrag}
\usepackage{soul}

\theoremstyle{definition}\newtheorem{definition}{Definition}[section]
\theoremstyle{definition}\newtheorem{proposition}[definition]{Proposition}
\theoremstyle{definition}\newtheorem{remark}[definition]{Remark}
\theoremstyle{definition}\newtheorem{theorem}[definition]{Theorem}
\theoremstyle{definition}
\theoremstyle{definition}
\theoremstyle{definition}\newtheorem{example}[definition]{Example}

\numberwithin{equation}{section}

\newcommand{\R}{\mathbb{R}}

\newcommand{\pp}[2]{\frac{\partial #1}{\partial #2}}

\newcommand{\n}{^{(n)}}
\newcommand{\ii}{^{(\infty )}}

\newcommand{\vv}{\mathbf{v}}

\newcommand{\pr}{\textbf{pr\,}}

\newcommand{\g}{\mathfrak{g}}

\newcommand{\D}{\mathcal{D}}

\newcommand{\G}{\mathcal{G}}

\newcommand{\J}{{\rm J}}

\newcommand{\oX}{\widetilde{X}}
\newcommand{\oY}{\widetilde{Y}}
\newcommand{\oU}{\widetilde{U}}
\newcommand{\overe}{\widetilde{e}}
\newcommand{\of}{\widetilde{f}}
\newcommand{\og}{\widetilde{g}}

\def\comp{\raise 1pt \hbox{$\,\scriptstyle\circ\,$}}

\begin{document}


\thispagestyle{fancy}
\fancyhead{}
\fancyfoot{}
\renewcommand{\headrulewidth}{0pt}
\cfoot{\thepage}
\rfoot{\today}

\vskip 1cm
\begin{center}
{\Large Invariant Discretization of Partial Differential Equations Admitting Infinite-Dimensional Symmetry Groups}
\vskip 1cm
Rapha\"el Rebelo$^{\dagger}$ and Francis Valiquette$^{\ddagger}$
\vskip 1cm
\end{center}

\noindent${}^{\dagger}$Centre de Recherche Math\'ematiques, Universit\'e de Montr\'eal, C.P. 6128, Succ. Centre-Ville, Montr\'eal, Qu\'ebec, H3C 3J7, Canada\\ 
Email: {\tt raph.rebelo@gmail.com}
\vskip 0.25cm
\noindent${}^{\ddagger}$Department of Mathematics, SUNY at New Paltz, New Paltz, NY 12561, USA\\
 Email:  {\tt valiquef@newpaltz.edu}

\vskip 0.5cm\noindent
{\bf Mathematics subject classification (MSC2010)}:  58J70, 65N06

\vskip 1cm

\abstract{ 
This paper is concerned with the invariant discretization of differential equations admitting infinite-dimensional symmetry groups.  By way of  example, we first show that there are differential equations with infinite-dimensional symmetry groups that do not admit enough joint invariants preventing the construction of invariant finite difference approximations.    To solve this shortage of joint invariants we  propose to discretize the pseudo-group action.  Computer simulations indicate that the numerical schemes constructed from the joint invariants of discretized pseudo-group can produce better numerical results than standard schemes.


\section{Introduction}

For the last 20 years, a considerable amount of work has been invested into the problem of invariantly discretizing differential equations with symmetries.   This effort is part of a larger program aiming to extend Lie's theory of transformation groups to finite difference equations, \cite{LW-2006}.  With the emergence of physical models based on discrete spacetime, and in light of the importance of symmetry in our understanding of modern physics, the problem of invariantly discretizing differential equations is still of present interest.  From a theoretical standpoint, working with invariant numerical schemes allows one to use standard Lie group techniques to find explicit solutions, \cite{RW-2004}, or compute conservation laws, \cite{DKW-2003}.  From a more practical point of view, the motivation stems from the fact that invariant schemes have been shown to outperform standard numerical methods in a number of examples, \cite{BRW-2008,D-2008,KO-2004,RV-2013}.

In general, to build an invariant numerical scheme one has to construct joint invariants (also known as finite difference invariants).  These joint invariants are usually found using one of two methods.  One can either use Lie's method of infinitesimal generators which requires solving a system of linear partial differential equations, \cite{D-2011,LW-2006}, or the method of equivariant moving frames which requires solving a system of (nonlinear) algebraic equations, \cite{KO-2004,O-2001}.   Both approaches produce joint invariants which, in the coalescent limit, converge to differential invariants of the prolonged action.   Thus far, the theory and applications found in the literature primarily deal with finite-dimensional Lie group actions and the case of infinite-dimensional Lie pseudo-groups as yet to be satisfactorily treated.  Many partial differential equations in hydrodynamics or meteorology admit infinite-dimensional symmetry groups.  The Navier--Stokes equation, \cite{O-1993}, the Kadomtsev--Petviashvili equation, \cite{DKLW-1985}, and the Davey--Stewartson equations, \cite{CW-1987}, are classical examples of such equations.  Linear or linearizable partial differential equations also form a large class of equations admitting infinite-dimensional symmetry groups.

To construct invariant numerical schemes of differential equations admitting symmetries, one of the main steps consists of finding joint invariants that approximate the differential invariants of the symmetry group.  For finite-dimensional Lie group actions, this can always be done by considering the product action on sufficiently many points.   Unfortunately, as the next example shows, the same is not true for infinite-dimensional Lie pseudo-group actions.

\begin{example}\label{example intro}
Let $f(x) \in \D(\mathbb{R})$ be a local diffeormorphism of $\mathbb{R}$.  Throughout the paper we will use the infinite-dimensional pseudo-group
\begin{equation}\label{main pseudo-group}
X=f(x),\qquad Y=y,\qquad U=\frac{u}{f^\prime(x)},
\end{equation}
acting on $\mathbb{R}^3 \setminus \{ u=0 \}$, to illustrate the theory and constructions. The pseudo-group \eqref{main pseudo-group} was introduced by Lie, \cite[p.373]{L-1895}, in his study of second order partial differential equations integrable by the method of Darboux.  It also appears in Vessiot's work on group splitting and automorphic systems, \cite{V-1904}, in Kumpera's investigation of Lie's theory of differential invariants based on Spencer's cohomology, \cite{K-1975}, and recently in \cite{OP-2005,OP-2008,P-2008} to illustrate a new theoretical foundation of moving frames.   

The differential invariants of the pseudo-group action \eqref{main pseudo-group} can be found in \cite{OP-2008}.  One of these invariants is 
\begin{equation}\label{I}
I_{1,1}=\frac{u\, u_{xy} - u_x\, u_y}{u^3}.
\end{equation}
With \eqref{I} it is possible to form the partial differential equation
\begin{equation}\label{main pde}
\frac{u\, u_{xy} - u_x\, u_y}{u^3} = 1,
\end{equation}
which was used in \cite{P-2008} to illustrate the method of symmetry reduction of exterior differential systems.  

By construction, Equation \eqref{main pde} is invariant under the pseudo-group\footnote{Equation \eqref{main pde} admits a larger symmetry group given by $X=f(x)$, $Y=g(y)$, $U=u/(f^\prime(x)\, g^\prime(y))$, with $f$, $g \in \D(\R)$.  This pseudo-group is considered in Example \ref{joint moving frame larger pseudo-group}.} \eqref{main pseudo-group}.    To obtain an invariant discretization of \eqref{main pde}, an invariant approximation of the differential invariant \eqref{I} must be found.  To discretize the invariant \eqref{I}, the multi-index $(m,n) \in \mathbb{Z}^2$ is introduced to label sample points:
\begin{equation}\label{sampling points}
x_{m,n},\qquad y_{m,n},\qquad u_{m,n}=u(x_{m,n},y_{m,n}).
\end{equation}
Following the general philosophy, \cite{D-2011,KO-2004,LW-2006,O-2001}, the pseudo-group \eqref{main pseudo-group} induces the product action
\begin{equation}\label{product action example 2}
X_{m,n} = f(x_{m,n}),\qquad Y_{m,n}=y_{m,n},\qquad U_{m,n}=\frac{u_{m,n}}{f^\prime(x_{m,n})}
\end{equation}
on the discrete points \eqref{sampling points}.   On an arbitrary finite set of points, we claim that the only joint invariants are
\begin{equation}\label{y invariants}
Y_{m,n}=y_{m,n}.
\end{equation}
To see this, let $\mathcal{N}$  be a finite subset of $\mathbb{Z}^2$, and  assume $x_{m,n}\in \text{dom}\, f$ for $(m,n) \in \mathcal{N}$.  Since the components $x_{m,n}$ are generically distinct and $f \in \D(\mathbb{R})$ is an arbitrary local diffeomorphism, the \emph{pseudo-group parameters}
\begin{equation}\label{f-fprime}
f(x_{m,n})\quad\text{and}\quad f^\prime(x_{m,n})\quad\text{with}\quad (m,n)\in \mathcal{N}
\end{equation}
are independent.  Hence, as shown in \cite{IOV-2011}, the pseudo-group \eqref{product action example 2} shares the same invariants as its Lie completion 
\begin{equation}\label{completion}
X_{m,n} = f_{m,n}(x_{m,n}),\qquad Y_{m,n}=y_{m,n},\qquad U_{m,n}=\frac{u_{m,n}}{f_{m,n}^\prime(x_{m,n})},
\end{equation}
%
where for each different subscript $(m,n) \in \mathcal{N}$, the functions $f_{m,n} \in \D(\R)$ are functionally independent local diffeomorphisms\footnote{It is customary to use the notation $f_{m,n} = f(x_{m,n})$ to denote the value of the function $f(x)$ at the point $x_{m,n}$, and this is the convention used in Sections \ref{discrete section}, \ref{numerical scheme section}, and \ref{numerics section}. In equation \eqref{completion}, the subscript attached to the diffeomorphism $f_{m,n}(x_{m,n})$ has a different meaning.  Here, the subscript $(m,n)$ is used to denote different diffeomorphisms. Thus, the pseudo-group \eqref{product action example 2} is contained in the Lie completion \eqref{completion}.  This particular use of the subscript only occurs in \eqref{completion}.}.  For the Lie completion \eqref{completion}, it is clear that \eqref{y invariants} are the only admissible invariants.  Hence, generically, we conclude that it is not possible to approximate the differential invariant \eqref{I} by joint invariants.  

To construct additional joint invariants, invariant constraints on the independent variables $x_{m,n}$ need to be imposed to reduce the number of pseudo-group parameters \eqref{f-fprime}.  To reduce this number as much as possible, we assume that 
\begin{equation}\label{x constraint}
x_{m,n+1}=x_{m,n}.
\end{equation}
Equation \eqref{x constraint} is seen to be invariant under the product action \eqref{product action example 2} since
$$
X_{m,n+1}=f(x_{m,n+1})=f(x_{m,n})=X_{m,n}
$$
when \eqref{x constraint} holds.  Equation \eqref{x constraint} implies that $x_{m,n}$ is independent of the index $n$:
$$
x_{m,n}=x_m.
$$
To cover (a region of) the $xy$-plane,
$$
\Delta x_m=x_{m+1}-x_{m} \neq 0\qquad \text{and}\qquad \delta y_{m,n}=y_{m,n+1}-y_{m,n}\neq 0  
$$
must hold. Since the variables $y_{m,n}$ are invariant under the product action \eqref{product action example 2} we can, for simplicity, set
\begin{equation}\label{yn}
y_{m,n}=y_n = k\, n + y_0,
\end{equation}
where $k>0$ and $y_0$ are constants.  To respect the product action \eqref{product action example 2} we cannot require the step size $\Delta x_m = x_{m+1}-x_m$ to be constant as this is not an invariant assumption of the pseudo-group action.  Thus, in general, the mesh in the independent variables $(x_m, y_n)$ will be rectangular with variable step sizes in $x$, see Figure \ref{mesh}.

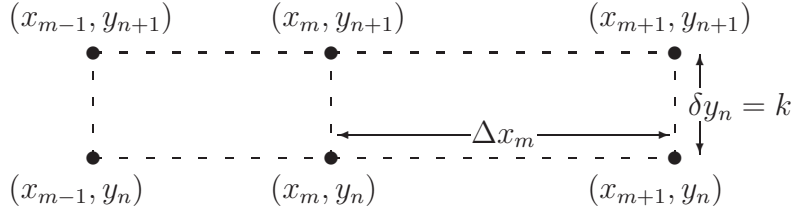
\begin{figure}[!h]
\begin{center}
\begin{picture}(280,80)
\put(30,20){\circle*{5}}
\put(-2,5){$(x_{m-1},y_n)$}
\put(120,20){\circle*{5}}
\put(97,5){$(x_m,y_n)$}
\put(250,20){\circle*{5}}
\put(217,5){$(x_{m+1},y_n)$}
\put(30, 60){\circle*{5}}
\put(-2,70){$(x_{m-1},y_{n+1})$}
\put(120,60){\circle*{5}}
\put(97,70){$(x_m,y_{n+1})$}
\put(250,60){\circle*{5}}
\put(217,70){$(x_{m+1},y_{n+1})$}
\multiput(30,20)(10,0){22}  {\line(1,0){3}}
\multiput(30,60)(10,0){22}  {\line(1,0){3}}
\multiput(30,20)(0,10){4}  {\line(0,1){3}}
\multiput(120,20)(0,10){4}  {\line(0,1){3}}
\multiput(250,20)(0,10){4}  {\line(0,1){3}}
\put(255,36.5){$\delta y_n=k$}
\put(259.5,47){\vector(0,1){12.5}}
\put(259.5,34){\vector(0,-1){12.5}}
\put(173,26){$\Delta x_m$}
\put(198,29){\vector(1,0){50}}
\put(172.5,29){\vector(-1,0){50}}
\end{picture}
\caption{Rectangular mesh.}\label{mesh}
\end{center}
\end{figure}  

Repeating the argument above, when \eqref{x constraint} and \eqref{yn} hold, the joint invariants of the product action \eqref{product action example 2} are
\begin{equation}\label{yu joint invariants}
y_{n},\qquad \frac{u_{m,n+k}}{u_{m,n}},\qquad k,\, m,\, n\in \mathbb{Z}.
\end{equation}
Introducing the dilation group
\begin{equation}\label{dilation}
X=x,\qquad Y=y,\qquad U=\lambda\, u,\qquad \lambda >0,
\end{equation}
we see that the differential invariant \eqref{I} cannot be approximated by the joint invariants \eqref{yu joint invariants}.  Indeed, since the 
 invariants $u_{m,n+k}/u_{m,n}$ are homogeneous of degree 0, any combination of the invariants \eqref{yu joint invariants} will converge to a differential invariant of homogeneous degree 0.  On the other hand, the differential invariant \eqref{I}  is homogeneous of degree $-1$ under  \eqref{dilation}.
\end{example}

As it stands, it is not possible to construct joint invariants that approximate the differential invariant \eqref{I}.  To remedy the problem, one possibility is to reduce the size of the symmetry group by considering sub-pseudo-groups.  For the diffeomorphism pseudo-group $\D(\mathbb{R})$, since the largest non-trivial sub-pseudo-group is the special linear group $SL(2)$, \cite{O-1995}, this approach drastically changes the nature of the action as it transitions from an infinite-dimensional transformation group to a three-dimensional group of transformations.  In this paper we are interested in preserving the infinite-dimensional nature of transformation groups and propound another suggestion.   Taking the point of view that the notion of derivative is not defined in the discrete setting, we propose to discretize infinite-dimensional pseudo-group actions.  In other words, derivatives are to be replaced by finite difference approximations. For the pseudo-group \eqref{main pseudo-group}, instead of considering the product action \eqref{product action example 2}, we suggest to work with the first order approximation
\begin{equation}\label{discrete pseudo-group example}
X_{m}=f(x_m),\qquad Y_{n}=y_{n}=k\, n + y_0,\qquad U_{m,n}=u_{m,n}\cdot \frac{x_{m+1}-x_m}{f(x_{m+1}) - f(x_m)}.
\end{equation}
In Section \ref{discrete section}, joint invariants of the pseudo-group action \eqref{discrete pseudo-group example} are constructed and an invariant numerical scheme approximating \eqref{main pde} is obtained in Section \ref{numerical scheme section}.

To develop our ideas we opted to use the theory of equivariant moving frames, \cite{O-2001,OP-2008}, but our constructions can also be recast within Lie's infinitesimal framework.  In Section \ref{pseudo-group section}, the concept of an infinite-dimensional Lie pseudo-group is recalled and the equivariant moving frame construction is summarized.  In Section \ref{discrete section}, pseudo-group actions are discretized and the equivariant moving frame construction is adapted to those actions.  Along with \eqref{main pseudo-group}, the pseudo-group
\begin{equation}\label{fg pseudo-group}
X=f(x),\qquad Y=y\, f^\prime(x) + g(x),\qquad U=u+\frac{y\, f^{\prime\prime}(x)+g^\prime(x)}{f^\prime(x)},
\end{equation}
with $f \in \D(\R)$ and $g \in C^\infty(\R)$, will stand as a second example to illustrate our constructions.   We choose to work with the pseudo-groups \eqref{main pseudo-group} and \eqref{fg pseudo-group} to keep our examples relatively simple.  Furthermore, these pseudo-groups have been extensively used in \cite{OP-2005,OP-2008,OP-2009-1,OP-2009-2} to illustrate the (continuous) method of moving frames.  With these well-documented examples, it allowed us to verify that our discrete constructions and computations did  converge to their continuous counterparts.

Finally, in Section \ref{numerics section} an invariant numerical approximation of  \eqref{main pde} is compared to a standard discretization of the equation. Our numerical tests show that the invariant scheme is more precise and stable than the standard scheme.  



\section{Lie Pseudo-groups and moving frames}\label{pseudo-group section}

For completeness, we begin by recalling the definition of a pseudo-group,  \cite{C-1953,K-1975,K-1959,OP-2005,OP-2008,SS-1965}. Let $M$ be an $m$-dimensional manifold.  By a local diffeomorphism of $M$ we mean a one-to-one map $\varphi\colon U \to V$ defined on open subsets $U$, $V=\varphi(U) \subset M$, with inverse $\varphi^{-1}\colon V \to U$.  

\begin{definition}\label{definition pg}
A collection  $\G$ of local  diffeomorphisms of $M$ is a \emph{pseudo-group} if
\begin{itemize}
\itemsep=-2pt
\item 
$\G$ is closed under restriction:  if $U\subset M$ is an open set and $g \colon  U \to M$ is in $\G$, then so is the restriction $g|_V$ for all open $V\subset U$.
\item
Elements of $\G$ can be pieced together:  if $U_\nu\subset M$ are open subsets, $U = \bigcup_\nu\, U_\nu$, and $g\colon U\to M$ is a local diffeomorphism with $g|_{U_\nu} \in \G$ for all $\nu $, then $g \in \G$.
\item
$\G$ contains the identity diffeomorphism: \, $\mathds{1} \cdot z = z$ for all $z \in M = \text{dom }\mathds{1}$.
\item
$\G$ is closed under composition: if $g \colon U\to M$ and $h\colon V \to M$ are two diffeomorphisms belonging to $\G$, and $g (U)\subset V$, then $h \cdot g \in \G$. 
\item
$\G$ is closed under inversion: if $g \colon U\to M$ is in $\G$ then so is $g^{-1}\colon g(U)\to M$.
\end{itemize}
\end{definition}

\begin{example}
One of the simplest pseudo-group is given by the collection of local diffeomorphisms $\D = \D(M)$ of a manifold $M$.  All other pseudo-groups defined on $M$ are sub-pseudo-groups of $\D$.  
\end{example}

For $0\leq n \leq \infty$, let $\D\n = \J^{(n)} (M,M)$ denote the bundle formed by the $n^\text{th}$ order jets of local diffeomorphisms of $M$.  Local coordinates on $\D\n$ are given by $\varphi^{(n)}|_z=(z,Z\n)$, where $z=(z^1,\ldots,z^m)$ are the source coordinates of the local diffeomorphism, $Z=\varphi(z)$, $Z=(Z^1,\ldots,Z^m)$ its target coordinates, and $Z\n$ collects the derivatives of the target coordinates $Z^a$ with respect to the source coordinates $z^b$ of order $\leq n$.  For $k\geq n$, the standard projection is denoted $\widetilde{\pi}^k_n\colon \D^{(k)}\to \D\n$.

\begin{definition}\label{Lie pseudo-group definition}
A pseudo-group $\G \subset \D$ is called a \emph{Lie pseudo-group} of order $n^\star \geq 1$ if, for all finite $n \geq n^\star\colon$
\begin{itemize}
\item $\G\n \subset \D\n$ forms a smooth embedded subbundle,
\item the projection $\widetilde{\pi}^{n+1}_n\colon \G^{(n+1)} \to \G\n$ is a fibration,
\item every local diffeomorphism $g \in \D$ satisfying $g^{(n^\star)} \subset \G^{(n^\star)}$ belongs to $\G$,
\item $\G\n = \text{pr}^{(n-n^\star)} \G^{(n^\star)}$ is obtained by prolongation.
\end{itemize}
\end{definition}

In local coordinates, the subbundle $\G^{(n^\star)} \subset \D^{(n^\star)}$ is characterized by a system of $n^{\star\, \text{th}}$ order (formally integrable) partial differential equations
\begin{equation}\label{determining system}
F^{(n^\star)}(z,Z^{(n^\star)})=0,
\end{equation}
called the $n^{\star\,\text{th}}$ order \emph{determining system} of the pseudo-group.  A Lie pseudo-group is said to be of \emph{finite type} if the solution space of \eqref{determining system} only involves a finite number of arbitrary constants.  Lie pseudo-groups of finite type are thus isomorphic to local Lie group actions.  On the other hand, a Lie pseudo-group is of \emph{infinite type} if it involves arbitrary functions.

\begin{remark}
Linearizing \eqref{determining system} at the identity jet $\mathds{1}^{(n^\star)}$ yields the \emph{infinitesimal determining equations}
\begin{equation}\label{infinitesimal determining system}
L^{(n^\star)}(z,\zeta^{(n^\star)})=0
\end{equation}
for an infinitesimal generator 
\begin{equation}\label{vector field}
\vv=\sum_{a=1}^m\> \zeta^a(z) \pp{}{z^a}.
\end{equation}
The vector field \eqref{vector field} is in the Lie algebra $\g$ of infinitesimal generators of $\G$ if its components are solution of \eqref{infinitesimal determining system}.  Given a differential equation $\Delta(x,u\n)=0$ with symmetry group $\G$, the infinitesimal determining system \eqref{infinitesimal determining system} is equivalent to the equations obtained by Lie's standard algorithm for determining the symmetry algebra of the differential equation $\Delta = 0$, \cite{O-1993}. 
\end{remark}

\begin{example}  
The pseudo-group \eqref{main pseudo-group} is a Lie pseudo-group.  The first order determining equations are
\begin{equation}\label{order 1 determining system}
\begin{gathered}
X_y=X_u=0,\qquad Y=y,\qquad Y_x=Y_u=0,\qquad Y_y=1,\\
U X_x = u,\qquad U_uX_x=1,\qquad U_y=0.
\end{gathered}
\end{equation}
If
$$
\vv = \xi(x,y,u)\pp{}{x} + \eta(x,y,u)\pp{}{y} + \varphi(x,y,u)\pp{}{u}
$$
denotes a (local) vector fields in $\R^3\setminus \{u=0\}$, the linearization of \eqref{order 1 determining system} at the identity jet yields the first order infinitesimal determining equations
$$
\xi_y=\xi_u=0,\quad \eta=\eta_x=\eta_y=\eta_u=0,\quad \varphi = -u\, \xi_x,\quad \varphi_y=0,\quad \varphi_u = - \xi_x.
$$
The general solution to this system of equations is
$$
\vv = a(x)\pp{}{x} - u\, a^\prime(x) \pp{}{u},
$$
where $a(x)$ is an arbitrary smooth function.
\end{example}

Given a Lie pseudo-group $\G$ acting on $M$, we are now interested in the induced action on $p$-dimensional submanifolds $S \subset M$ with $1 \leq p < m = \text{dim }M$.  It is customary to introduce adapted coordinates 
\begin{equation}\label{graph coordinates}
z=(x,u)=(x^1,\ldots,x^p,u^1,\ldots, u^q)
\end{equation}
on $M$ so that, locally, a submanifold $S$ transverse to the vertical fibre $\{x=x_0\}$ is given as the graph 
of a function $S=\{(x,u(x))\}$.   For each integer $0\leq n \leq \infty$, let $\J^{(n)}=\J^{(n)}(M,p)$ denote the $n^\text{th}$ order \emph{submanifold jet bundle} defined as the set of equivalence classes under the equivalence relation of $n^\text{th}$ order contact, \cite{O-1995}.  For $k\geq n$, let $\pi^k_n \colon \J^{(k)}\to \J^{(n)}$ denote the canonical projection.   In the adapted system of coordinates $z=(x,u)$, coordinates on $\J^{(n)}$ are given by
\begin{equation}\label{graph jet variables}
z\n=(x,u\n)=(\, \ldots\; x^i\, \ldots\, u^\alpha_{x^J}\, \ldots\,),
\end{equation}
where $u\n$ denotes the collection of derivatives $u^\alpha_{x^J}$ of order $0 \leq \# J \leq n$. 

Alternatively, when no distinction between dependent and independent variables is made, a submanifold $S \subset M$ can be locally parameterized by $p$ variables $s=(s^1,\ldots,s^p)\in \mathbb{R}^p$ so that
$$
z(s)=(x(s),u(s))\in S.
$$
In the numerical analysis community, the variables $s=(s^1,\ldots,s^p)$ are called \emph{computational variables}, \cite{HR-2011}.  We let $\mathcal{J}\n$ denote the $n^\text{th}$ order jet space of submanifolds $S\subset M$ parametrized by computational variables.  Local coordinates on $\mathcal{J}\n$ are given by
\begin{equation}\label{computational jet coordinates}
\mathfrak{z}\n=(s,x\n,u\n)=(\, \ldots\, s^i\, \ldots\, x^i_{s^A}\, \ldots\, u^\alpha_{s^A}\, \ldots\,),
\end{equation}
with $1\leq i \leq p$, $1\leq \alpha \leq q$, and $0 \leq \# A \leq n$.  The transition between the jet coordinates \eqref{graph jet variables} and \eqref{computational jet coordinates} is given by the chain rule.  Provided
\begin{equation}\label{change of variables condition}
\text{det}\bigg( \pp{x^i}{s^j} \bigg) \neq 0,
\end{equation}
successive application of the implicit total differential operators
\begin{equation}\label{implicit operators}
D_{x^i} = \sum_{j=1}^p\> W^j_i\, D_{s^j},\qquad \big(W^j_i\big)=\big(x_{s^i}^j\big)^{-1},
\end{equation}
to the dependent variables $u^\alpha$ will give the coordinate expressions for the $x$ derivatives of $u$ in terms of the $s$ derivatives of $x$ and $u$:
\begin{equation}\label{uxJ}
u^\alpha_{x^J} = D_{x^{j_1}}\cdots D_{x^{j_k}} u^\alpha = \bigg(\sum_{\ell=1}^p\> W^\ell_{j_1} \, D_{s^\ell}\bigg) \cdots \bigg(\sum_{\ell=1}^p\> W^\ell_{j_k}\, D_{s^\ell}\bigg) u^\alpha.
\end{equation}

Given a Lie pseudo-group $\G$ acting on $M$, the action is prolonged to the computational variables by requiring that they remain unchanged:
$$
g\cdot (s,z) = (s,g\cdot z)\qquad \text{for all}\qquad g \in \G.
$$
By abuse of notation we still use $\G$ to denote the extended action $\{ \mathds{1} \} \times \G$ on $\mathbb{R}^p \times M$.   

The complete theory of moving frames for infinite-dimensional Lie pseudo-groups can be found in \cite{OP-2008}.  For reasons that will become more apparent in the next section we recall the main constructions over the jet bundle $\mathcal{J}\n$ rather than $\J\n$.  Using \eqref{uxJ} one can translate the constructions from $\mathcal{J}\n$ to $\J\n$.  Let
\begin{equation}\label{lifted bundle}
\mathcal{B}\n = \mathcal{J}\n \times_{M} \G\n
\end{equation}
denote the $n^{\text{th}}$ order \emph{lifted bundle}.  Local coordinates on $\mathcal{B}\n$ are given by $(\mathfrak{z}\n,g\n)$, where the base coordinates are the submanifold jet coordinates $\mathfrak{z}\n=(s,x\n,u\n)\in \mathcal{J}\n$ and the fibre coordinates are the pseudo-group parameters $g\n$ where $(x,u) \in \text{dom }g$.  A local diffeomorphism $h \in \G$  acts on $\mathcal{B}\n$ by right multiplication:
\begin{equation}\label{lifted action}
R_h (\mathfrak{z}\n,g\n)=(h\n \cdot \mathfrak{z}\n, g\n \cdot (h\n)^{-1}),
\end{equation}
where defined.  The second component of  \eqref{lifted action} corresponds to the usual right multiplication $R_h(g\n) = g\n \cdot (h\n)^{-1}$ of the pseudo-group onto $\G\n$, \cite{OP-2008}.  The first component $h\n \cdot \mathfrak{z}\n = (s,h\n \cdot x\n, h\n\cdot u\n) = (s,X\n,U\n)$ is the prolonged action of the pseudo-group $\G$ onto the jet space $\mathcal{J}\n$. Coordinate expressions for the prolonged action are obtained by differentiating the target coordinates $Z=(X,U)$ with respect to the computational variables $s$:
\begin{equation}\label{prolonged action}
X^i_A=D^{a^1}_{s^1}\ldots D^{a^p}_{s^p}  X^i,\qquad U^\alpha_A = D^{a^1}_{s^1}\ldots D^{a^p}_{s^p} U^\alpha,
\end{equation}
where $A=(a^1,\ldots,a^p)$. The expressions \eqref{prolonged action} are invariant under the \emph{lifted action} \eqref{lifted action} and these functions are called \emph{lifted invariants}.

\begin{definition}
A \emph{(right) moving frame} of order $n$ is a $\G$-equivariant section $\widehat{\rho}^{\,(n)}$ of the lifted bundle $\mathcal{B}\n \to \mathcal{J}\n$. 
\end{definition}

In local coordinates, the notation 
$$
\widehat{\rho}^{\,(n)}(\mathfrak{z}\n)=(\mathfrak{z}\n, \rho\n(\mathfrak{z}\n))
$$
is used to denote an order $n$ right moving frame.  Right equivariance means that for $g \in \G$
$$
R_g \widehat{\rho}^{\,(n)}(\mathfrak{z}\n) = \widehat{\rho}^{\,(n)}(g\n\cdot \mathfrak{z}\n),
$$
where defined.

\begin{definition}
Let 
$$
\G\n_{\mathfrak{z}\n}= \bigg\{ g\n \in \G\n|_{\mathfrak{z}} :\, g\n\cdot \mathfrak{z}\n = \mathfrak{z}\n \bigg\}
$$
denote the \emph{isotropy subgroup} of $\mathfrak{z}\n \in \mathcal{J}\n$.  The pseudo-group $\G$ is said to act \emph{freely} at $\mathfrak{z}\n \in \mathcal{J}\n$ if $\G\n_{\mathfrak{z}\n}=\{ \mathds{1}\n|_\mathfrak{z}\}$.  The pseudo-group $\G$ is said to act \emph{freely at order $n$} if it acts freely on an open subset $\mathcal{V}^{(n)} \subset \mathcal{J}\n$, called the set of \emph{regular $n$-jets.}
\end{definition}

\begin{theorem}\label{existence moving frame}
Suppose $\G$ acts freely on $\mathcal{V}^{(n)} \subset \mathcal{J}\n$, with its orbits forming a regular foliation.
Then an $n^\text{th}$ order moving frame exists in a neighbourhood of $\mathfrak{z}\n \in \mathcal{V}^{(n)}$.
\end{theorem}

Once a pseudo-group acts freely, a result known as the \emph{persistence of freeness}, \cite{OP-2009-1,OP-2009-2}, guarantees that the action remains free under prolongation. 

\begin{theorem}\label{persistence of freeness}
If a Lie pseudo-group $\G$ acts freely at $\mathfrak{z}\n$, then it acts freely at any $\mathfrak{z}^{(k)} \in \mathcal{J}^{(k)}$, $k\geq n$, with $\pi^k_n(\mathfrak{z}^{(k)}) = \mathfrak{z}\n$.  
\end{theorem}

\begin{remark}
Theorems \ref{existence moving frame} and \ref{persistence of freeness} also hold when the pseudo-group action is locally free, meaning that the isotropy group $\G\n_{\mathfrak{z}\n}$ is a discrete subgroup of $\G\n\big|_{\mathfrak{z}}$.   
\end{remark}

An order $n \geq n^\star$ moving frame is constructed through a normalization procedure based on the choice of a cross-section $\mathcal{K}^{(n)} \subset \mathcal{V}^{(n)}$ to the pseudo-group orbits.  The associated (locally defined) right moving frame section $\widehat{\rho}^{\,(n)}\colon \mathcal{V}^{(n)} \to \mathcal{B}\n$ is uniquely characterized by the condition that $\rho\n(\mathfrak{z}\n)\cdot \mathfrak{z}\n \in \mathcal{K}^{(n)}$.  In coordinates, assuming that 
\begin{equation}\label{coordinate cross-section}
\mathcal{K}^{(n)}=\{z_{i_1}=c_1, \ldots,z_{i_{r_n}}=c_{r_n}\,:\, r_n=\text{dim }\G\n|_{\mathfrak{z}}\}
\end{equation}
is a coordinate cross-section, the moving frame $\widehat{\rho}^{\,(n)}$ is obtained by solving the \emph{normalization equations}
$$
Z_{i_1}(s,x\n,u\n,g\n)=c_1,\qquad \ldots \qquad Z_{i_{r_n}}(s,x\n,u\n,g\n)=c_{r_n},
$$
for the pseudo-group parameters $g\n=\rho\n(\mathfrak{z}\n)$.    As one increases the order from $n$ to $k>n$, a new cross-section $\mathcal{K}^{(k)} \subset \J^{(k)}$ must be selected.   These cross-sections are required to be compatible meaning that $\pi^k_n(\mathcal{K}^{(k)}) = \mathcal{K}^{(n)}$ for all $k > n$.  This in turn, implies the compatibility of the moving frames: $\widehat{\pi}^k_n \comp \widehat{\rho}^{\,(k)} = \widehat{\rho}^{\,(n)}\comp \pi^k_n$, where $\widehat{\pi}^k_n\colon \mathcal{B}^{(k)} \to \mathcal{B}\n$ is the standard projection.

\begin{proposition}
Let $\widehat{\rho}^{\,(n)}$ be an $n^\text{th}$ order right moving frame.  The \emph{normalized invariants} 
$$
(s,H\n,I\n)=\iota(s,x\n,u\n)=\rho\n(\mathfrak{z}\n) \cdot \mathfrak{z}\n,
$$
form a complete set of differential invariants of order $\leq n$.
\end{proposition}

\begin{example}\label{continuous moving frame example}
In this example we construct a moving frame for the pseudo-group \eqref{main pseudo-group}.  The computations for graphs of functions $(x,y,u(x,y))$ appear in \cite{OP-2008}.  In preparation for the next section we revisit the calculations using the computational variables $(s,t)$ so that $x=x(s,t)$, $y=y(s,t)$ and $u=u(s,t)$.   To simplify the computations let
\begin{subequations}\label{continuous xy constraints}
\begin{equation}
y = k\, t + y_0,
\end{equation}
where $k > 0$ and $y_0$ are constants, and assume that
\begin{equation}\label{x-constraint}
x_t=0.
\end{equation}
\end{subequations}
In other words, $x=x(s)$ is a function of the computational variable $s$.  We note that the constraints \eqref{continuous xy constraints} are invariant under the pseudo-group action \eqref{main pseudo-group}. For the $y$ variable, this is straightforward as it is an invariant.  The invariance of \eqref{x-constraint} follows from the chain rule:
$$
X_t = f_x\, x_t = 0\qquad \text{when}\qquad x_t=0.
$$
The non-degeneracy condition \eqref{change of variables condition} requires the invariant constraint $x_s \neq 0$ to be satisfied.  

Up to order 2, the prolonged action is
\begin{equation}\label{prolonged action main pseudo-group}
\begin{gathered}
S=s,\qquad T=t,\qquad  Y=y,\qquad X=f(x),\qquad U=\frac{u}{f_x},\\
X_s = f_x\, x_s,\qquad Y_t=k,\qquad U_s = \frac{u_s}{f_x}-\frac{u\, f_{xx}\,x_s}{f_x^2},\qquad U_t = \frac{u_t}{f_x},\\
X_{ss} = f_{xx}\, x_s^2 + f_x\, x_{ss},\qquad Y_{tt}=0,\qquad  U_{tt} = \frac{u_{tt}}{u},\qquad U_{st} = \frac{u_{st}}{f_x} - \frac{u_t\, f_{xx}\, x_s}{f_x^2},\\
U_{ss} = \frac{u_{ss}}{f_x} + 2 \frac{u\, f_{xx}^2\, x_s^2}{f_x^3} - 2 \frac{u_s\, f_{xx}\, x_s}{f_x^2} - \frac{u\, f_{xxx}\, x_s^2}{f_x^2} - \frac{u\, f_{xx}\, x_{ss}}{f_x^2}.
\end{gathered}
\end{equation}
A cross-section to the prolonged action \eqref{prolonged action main pseudo-group} and its prolongation is given by
\begin{equation}\label{cross-section}
\mathcal{K}^{(\infty)} = \{x=0,\; u=1,\; u_{s^k}=0,\; k\geq 1\}.
\end{equation}
Solving the normalization equations
$$
X=0,\qquad U=1,\qquad U_{s^k}=0,\qquad k\geq 1,
$$
for the pseudo-group parameters $f, f_x, f_{xx}, \ldots,$ we obtain the right moving frame
$$
f=0,\qquad f_x = u,\qquad f_{xx} = \frac{u_s}{x_s},\qquad f_{xxx}=\frac{u_{ss}}{x_s^2} - \frac{u_s\, x_{ss}}{x_s^3},\qquad\ldots.
$$
In general,
\begin{equation}\label{main pseudo-group normalizations}
f=0,\qquad f_{x^{k+1}}=\bigg(\frac{D_s}{x_s}\bigg)^{k}u,\qquad k\geq 0.
\end{equation}
Substituting the pseudo-group normalizations \eqref{main pseudo-group normalizations} into the prolonged action \eqref{prolonged action main pseudo-group} yields the normalized differential invariants
\begin{equation}\label{normalized invariants pseudo-group1}
\begin{gathered}
\iota(s)=s,\qquad \iota(t)=t,\qquad \iota(y)=y,\qquad I_1=\iota(x_s) = u\, x_s,\qquad \iota(y_t)=k,\\
J_{0,1}=\iota(u_t)=\frac{u_t}{u},\qquad I_2=\iota(x_{ss})=u_s\, x_s + u\, x_{ss},\qquad \iota(y_{tt})=0,\\ 
J_{0,2}=\iota(u_{tt})=\frac{u_{tt}}{u},\qquad J_{1,1}=\iota(u_{st})=\frac{u\, u_{st}-u_t\, u_s}{u^2}.
\end{gathered}
\end{equation}
\end{example}

\begin{remark}
To transition between the expressions obtained in Example \ref{continuous moving frame example} and those appearing in \cite{OP-2008}, it suffices to use the chain rule.   When \eqref{continuous xy constraints} holds, 
$$
D_x = \frac{D_s}{x_s},\qquad D_y = \frac{D_t}{k},
$$
so that
\begin{equation}\label{chain rule}
u_x = \frac{u_s}{x_s},\quad\; u_y=\frac{u_t}{k},\quad\; u_{xx} = \frac{x_s\, u_{ss} - u_s\, x_{ss}}{x_s^3},\quad\; u_{xy}=\frac{u_{st}}{k\, x_s},\quad\; u_{yy}=\frac{u_{tt}}{k^2},\quad\; \ldots.
\end{equation}
The prolonged action on $u_x$, $u_y$, $\ldots$, can then be obtained by substituting \eqref{prolonged action main pseudo-group} into \eqref{chain rule}.  For example,
$$
U_X = \frac{U_s}{X_s} = \bigg(\frac{u_s}{f_x} - \frac{u\, f_{xx}\, x_s}{f_x^2}\bigg)\frac{1}{f_x\, x_s} = \frac{u_x}{f_x^2} - \frac{u\, f_{xx}}{f_x^3}.
$$
In the jet variables $z\ii = (x, y, u\ii) = (x, y, u, u_x, u_y, \ldots)$ a cross-section is given by, \cite{OP-2008},
\begin{equation}\label{OP cross-section}
\overline{\mathcal{K}}^{(\infty)} = \{x=0,\, u=1,\, u_{x^k}=0,\, k\geq 1\},
\end{equation}
and the corresponding moving frame is 
\begin{equation}\label{OP moving frame}
f=0,\qquad f_{x^{k+1}}=u_{x^k},\qquad k\geq 0.
\end{equation}
Expressing $u_{x^k}$ in terms of the derivatives $x_{s^k}$, $u_{s^k}$ using \eqref{chain rule}, one sees that \eqref{OP cross-section} and \eqref{OP moving frame} are equivalent to \eqref{cross-section} and \eqref{main pseudo-group normalizations} in the computational variable framework.  In the following, the cross-sections \eqref{cross-section} and \eqref{OP cross-section} (and the corresponding moving frames \eqref{main pseudo-group normalizations} and \eqref{OP moving frame}) are said to be \emph{equivalent}.  

Not all cross-sections are equivalent.  For example, instead of using the cross-section \eqref{cross-section}, it is also possible to choose the (non-minimal) cross-section
\begin{equation}\label{incompatible cross-section}
\widetilde{\mathcal{K}}^{(\infty)} = \{ x=0,\, x_s = 1,\, x_{s^{k+2}}=0,\, k\geq 0 \}.
\end{equation}
Since \eqref{OP cross-section} is not related to \eqref{incompatible cross-section} by the substitutions \eqref{chain rule}, the cross-section \eqref{incompatible cross-section} is said to be \emph{inequivalent} to \eqref{OP cross-section}.
\end{remark}

\begin{definition}
Let $\G\n$ be a Lie pseudo-group acting on $\J\n$ and $\mathcal{J}\n$.  A cross-section $\mathcal{K}\n \subset \mathcal{J}\n$ is said to be \emph{equivalent} with  the cross-section $\overline{\mathcal{K}}\n \subset \J\n$ if the defining equation \eqref{coordinate cross-section} of $\mathcal{K}\n$ are obtained from those of $\overline{\mathcal{K}}\n$ by expressing the submanifold jet $z\n$ in terms of $\mathfrak{z}\n$ using the relations \eqref{uxJ}.  
\end{definition}


\section{Discrete pseudo-groups and moving frames}\label{discrete section}

Let $M^{\times k}$ denote the $k$-fold Cartesian product of a manifold $M$.   Discrete points in $M$ are labelled using the multi-index notation
\begin{equation}\label{multi-index notation}
z_N = (x_N,u_N),\qquad N=(n^1,\ldots,n^p)\in \mathbb{Z}^p.
\end{equation}
The multi-index notation \eqref{multi-index notation} can be related to the continuous theory of Section \ref{pseudo-group section} in the following way.   The multi-index $N=(n^1,\ldots,n^p)\in \mathbb{Z}^p \subset \mathbb{R}^p$ can be thought as sampling the computational variables $s=(s^1,\ldots,s^p)\in \mathbb{R}^p$ on a unit hypercube grid.  Thus, the notation $z_N = z(N)$ can be understood as sampling a submanifold $S \subset M$ parameterized by $z(s)=(x(s),u(s))$ at the integer points $N \in \mathbb{Z}^p$. 

To mimic the continuous theory of moving frames in the finite difference setting, a discrete counterpart to the submanifold jet space $\mathcal{J}\n$ is introduced.

\begin{definition}
Let $M$ be a manifold with local coordinate system $z=(x,u)$.  The $k$-fold \emph{joint product} of $M$ is a subset of the $k$-fold Cartesian product $M^{\times k}$ given by
$$
M^{\diamond k}=\{(z_{N_1},\ldots,z_{N_k})\,| \, z_{N_i} \neq z_{N_j}\, \text{for all } i\neq j\} \subset M^{\times k}.
$$
\end{definition}
 
\begin{definition}
The $n^\text{th}$ order \emph{forward discrete jet} at the multi-index $N$ is the point
\begin{equation}\label{discrete jet}
\mathfrak{z}^{[n]}_N = (N, \ldots\, z_{N+K}\, \ldots),
\end{equation}
where $(\ldots\, z_{N+K}\, \ldots) \in M^{\diamond d_n}$ with 
$$
d_n=m \binom{p+n}{n},\qquad m = \text{dim }M,
$$
and $K=(k^1,\ldots , k^p)$ is a non-negative multi-index of order $0 \leq \# K \leq n$. 
\end{definition}

In dimension 2, when $N=(m,n)$, Figure \ref{forward joint spaces} shows the multi-indices contained in a forward discrete jet or order $\leq 2$.  In general, the multi-indices included in $\mathfrak{z}^{[k]}_{m,n}$ are those contained in the interior and boundary of the right isosceles triangle with vertices at $(m,n)$, $(m+k,n)$ and $(m,n+k)$.  

\begin{figure}[!h]
\centering
\subfloat[$k=0$]{\setlength{\unitlength}{0.02cm}
\begin{picture}(100,170)
\put(0,30){\vector(1,0){150}}
\put(10,20){\vector(0,1){150}}
\put(150,32){$s$}
\put(12,165){$t$}
\put(10,30){\circle*{10}}
\put(-12,12){{\scriptsize$(m,n)$}}
\end{picture}}
\hskip 2cm
\subfloat[$k=1$]{\setlength{\unitlength}{0.02cm}
\begin{picture}(100,170)
\put(0,30){\vector(1,0){150}}
\put(10,20){\vector(0,1){150}}
\put(150,32){$s$}
\put(12,165){$t$}
\put(10,30){\circle*{10}}
\put(10,92.5){\circle*{10}}
\put(72.5,30){\circle*{10}}
\put(-12,12){{\scriptsize$(m,n)$}}
\put(40,12){{\scriptsize$(m+1,n)$}}
\put(-62,90){{\scriptsize$(m,n+1)$}}
{\thicklines
\color{red}\put(10,90){\line(0,-1){60}}
\color{red}\put(10,30){\line(1,0){65}}
\color{red}\put(10,92){\line(1,-1){62.5}}}
\end{picture}}
\hskip 2cm
\subfloat[$k=2$]{\setlength{\unitlength}{0.02cm}
\begin{picture}(150,170)
\put(0,30){\vector(1,0){150}}
\put(10,20){\vector(0,1){150}}
\put(150,32){$s$}
\put(12,165){$t$}
\multiput(10,92.5)(10,0){7} {\line(1,0){5}}
\multiput(72.5,27.5)(0,10){7}{\line(0,1){5}}
\put(10,30){\circle*{10}}
\put(10,92.5){\circle*{10}}
\put(72.5,30){\circle*{10}}
\put(72.5,92.5){\circle*{10}}
\put(10,155){\circle*{10}}
\put(135,30){\circle*{10}}
\put(-12,12){{\scriptsize$(m,n)$}}
\put(40,12){{\scriptsize$(m+1,n)$}}
\put(-62,90){{\scriptsize$(m,n+1)$}}
\put(80,90){{\scriptsize$(m+1,n+1)$}}
\put(112.5,12){{\scriptsize$(m+2,n)$}}
\put(-62,150){{\scriptsize$(m,n+2)$}}
{\thicklines
\color{red}\put(78,88){\line(1,-1){58}}
\color{red}\put(78,88){\line(-1,1){67}}
\color{red}\put(10,155){\line(0,-1){125}}
\color{red}\put(10,30){\line(1,0){125}}}
\end{picture}}
\caption{Multi-indices occurring in $\mathfrak{z}^{[k]}_{m,n} \in \mathcal{J}^{[k]}$ for $k=0,1,2$.}\label{forward joint spaces}
\end{figure}
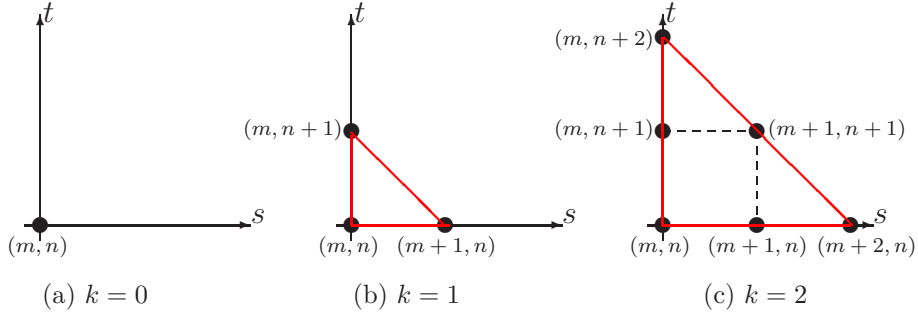

\begin{definition}
The $n^\text{th}$ order \emph{forward joint jet space} $\mathcal{J}^{[n]}$ is the collection of forward discrete jets \eqref{discrete jet}:
$$
\mathcal{J}^{[n]} = \bigcup_{N\,\in\, \mathbb{Z}^p}\, \mathfrak{z}^{[n]}_N.
$$
\end{definition}

For $k >n$, $\pi^k_n\colon \mathcal{J}^{[k]} \to \mathcal{J}^{[n]}$ will denote the projection obtained by truncating $\mathfrak{z}^{[k]}_N = (N, \ldots\,z_{N+K}\,\ldots)$, $0\leq \# K \leq k$, to
$$
\pi^k_n(\mathfrak{z}^{[k]}_N)=\mathfrak{z}^{[n]}_N= (N,\ldots\,z_{N+K}\, \ldots),\qquad 0\leq \# K \leq n.
$$

Let us explain how $\mathcal{J}^{[n]}$ can be understood as a discrete representation of the submanifold jet space $\mathcal{J}^{(n)}$.  For this, let $e_i=(0,\ldots,0,1,0,\ldots,0)$ be the $i^\text{th}$ element of the standard orthonormal basis of $\mathbb{R}^p$, and let
$$
S_i (N) = N+e_i,\qquad i=1,\ldots,p,
$$
denote the forward shift operator in the $i^\text{th}$ component.  Then, on a unit hypercube grid in the computational variables,  the derivative operators $D_{s^i}$ can be approximated by the forward difference
$$
D_{s^i} \sim \Delta_i = S_i - \mathds{1},\qquad i=1,\ldots p,
$$
where $\mathds{1}(N) = N$ is the identity map.  Then, for a non-negative multi-index $K=(k^1,\ldots,k^p)$,
\begin{equation}\label{z_s^k approximation}
z^N_{s^K} = \Delta_1^{k_1} \cdots \Delta_p^{k_p}(z_N)
\end{equation}
is a forward difference approximation of the derivative $z_{s^K}$ at the point $s=N$.  Making the change of variables $z_{N+K} \mapsto z_{s^K}^N$, we have that 
$$
\mathfrak{z}^{[n]}_N \simeq (N, \ldots\, z^N_{s^K}\, \ldots\,)=(N, \ldots\, x^N_{s^K} \ldots u^N_{s^K}\ldots\,),\qquad 0\leq \# K \leq n,
$$
is a finite difference approximation of the submanifold jet $(s,x^{(n)},u^{(n)})$ at the point $s=N$ on a unit hypercube grid.  In this sense, $\mathfrak{z}^{[n]}_N$ can be thought as a discrete counterpart to the submanifold jet $\mathfrak{z}\n=(N,x^{(n)},u^{(n)})$ in the computational variable formalism, \cite{HR-2011}.   

\begin{remark}
In \eqref{z_s^k approximation} and elsewhere, the usual derivative notation is supplemented by a superscript to denote (forward) discrete derivatives.  The superscript indicates where the derivative is evaluated.
\end{remark}

\begin{remark}
It is also possible to introduce a backward discrete jet space by introducing the backward differences
$$
\nabla_i = \mathds{1} - S_i^{-},\qquad\text{where}\qquad S_i^{-}(N) = N - e_i.
$$
For numerical purposes, it might be preferable to consider symmetric discrete jets, but to simplify the exposition we restrict ourself to forward differences.  All constructions can be adapted to these alternative  settings.
\end{remark}

Now, assume that the discrete counterpart of the non-degeneracy condition \eqref{change of variables condition} holds. Namely,
\begin{equation}\label{mesh non-degeneracy condition}
\text{det }(\Delta_j(x^i_N)) \neq 0.
\end{equation}
Then, discrete approximations $u^{\alpha;N}_{x^J}$ of the derivatives $u^\alpha_{x^J}$ can be obtained as follows:
\begin{enumerate}
\item compute the expressions \eqref{uxJ},
\item replace the derivatives $D_{s^i}$ by the difference operators $\Delta_i$.
\end{enumerate}
Since the independent variables $x^i_N$ do not have to form a rectangular grid, the finite difference approximations $u^{\alpha;N}_{x^J}$ will hold on any admissible mesh.   Having these expressions will be important as below a Lie pseudo-group will act on $z_N=(x_N,u_N)$ and the expressions for $u_{x^J}^{\alpha;N}$ need to hold on general meshes, \cite{D-2011,LW-2006,RV-2013}.   

%
%

Using the approximations $u_{x^J}^{\alpha;N}$, a finite difference approximation of the jet space $\J\n$ is given by
$$
\J\n \sim \J^{[n]} = \bigcup_{N\in \mathbb{Z}^p} (x_N,\ldots\, u^{\alpha;N}_{x^J}\, \ldots\,),\qquad 0\leq \# J \leq n.
$$

\begin{example}
To illustrate the above discussion, we consider the case of two independent variables $(x,y)$ and one dependent variable $u(x,y)$.  Introducing the computational variables $(s,t)\in \mathbb{R}^2$ so that $x=x(s,t)$ and $y=y(s,t)$, the implicit total derivative operators \eqref{implicit operators} are
\begin{equation}\label{xy implicit differentiation}
D_x = \frac{y_t\, D_s - y_s\, D_t}{x_s\, y_t - y_s\, x_t},\qquad
D_y = \frac{x_s\, D_t-x_t\, D_s}{x_s\, y_t - y_s\, x_t},
\end{equation}
with $x_s\, y_t - y_s\, x_t \neq 0$.  Applying \eqref{xy implicit differentiation} to the dependent variable $u$ yields
\begin{equation}\label{first order derivatives}
u_x = D_x\, u = \frac{y_t\, u_s - y_s\, u_t}{x_s\, y_t - y_s\, x_t}\qquad \text{and}\qquad u_y = D_y\, u = \frac{x_s\, u_t - x_t\, u_s}{x_s\, y_t - y_s\, x_t}.
\end{equation}
Using the multi-index $N=(m,n)\in \mathbb{Z}^2 \subset \mathbb{R}^2$  to sample the computational variables $(s,t)$ at integer values and introducing the shift operators
$$
S_1 (m,n) = (m+1,n),\qquad S_2 (m,n)=(m,n+1),
$$
and the difference operators
\begin{equation}\label{2D difference operators}
D_s \sim \Delta = S_1 - \mathds{1},\qquad D_t \sim \delta = S_2 - \mathds{1},
\end{equation}
finite difference approximations of the first order partial derivatives \eqref{first order derivatives} are given by 
\begin{equation}\label{first order discrete derivatives}
u_x^{m,n} = \frac{\delta y_{m,n}\, \Delta u_{m,n} - \Delta y_{m,n}\, \delta u_{m,n}}{\Delta x_{m,n}\, \delta y_{m,n} - \Delta y_{m,n}\, \delta x_{m,n}},\qquad
u_y^{m,n} = \frac{\Delta x_{m,n}\, \delta u_{m,n} - \delta x_{m,n}\, \Delta u_{m,n}}{\Delta x_{m,n}\, \delta y_{m,n} - \Delta y_{m,n} \delta x_{m,n}},
\end{equation}
provided $\Delta x_{m,n}\, \delta y_{m,n} - \Delta y_{m,n}\, \delta x_{m,n} \neq 0$.

The expressions \eqref{first order derivatives} and their finite difference approximations \eqref{first order discrete derivatives} can be simplified if constraints on the functions $x(s,t)$ and $y(s,t)$ are imposed.  For example, in Example \ref{pseudo-group discretization 2} we will impose the constraints
\begin{equation}\label{xy differential constraints}
x_t=0\qquad \text{and}\qquad y_{tt}=0,
\end{equation}
so that $x=x(s)$ and $y=t\,f(s)+g(s)$, with $f(s)\cdot x^\prime(s)\neq 0$. The operators \eqref{xy implicit differentiation} then reduce to
\begin{equation}\label{DxDy - DsDt}
D_x = \frac{y_t\ D_s - y_s\, D_t}{x_s\, y_t},\qquad D_y = \frac{D_t}{y_t},
\end{equation}
and
\begin{gather}
u_x = \frac{y_t\,u_s-y_s\, u_t}{x_s\, y_t},\qquad u_y=\frac{u_t}{y_t},\qquad 
u_{yy}=\frac{u_{tt}}{y_t^2},\qquad
u_{xy}=\frac{y_t\, u_{st}-y_{st}\, u_t - y_s\, u_{tt}}{x_s\, y_t^2},\nonumber\\
u_{xx}=\frac{1}{x_s}\bigg[\frac{y_{st}\, u_s+y_t\, u_{ss} - y_{ss}\, u_t - y_s\, u_{st}-(x_{ss}\,y_t + x_s\, y_{st})u_x}{x_s\, y_t} - y_s\, u_{xy}\bigg],\label{derivative computational variable expressions}\\
u_{yyy}=\frac{u_{ttt}}{y_t^3},\qquad
u_{xyy}=\frac{y_t\, u_{stt} - 2 u_{tt}\, y_{st} - y_s\, u_{ttt}}{x_s\, y_t^3},\qquad \ldots.\nonumber
\end{gather}
At the discrete level the differential constraints \eqref{xy differential constraints} are replaced by
$$
\delta x_{m,n} = 0\qquad \text{and}\qquad \delta^2 y_{m,n} = \delta y_{m,n+1} - \delta y_{m,n} = y_{m,n+2} - 2 y_{m,n+1} + y_{m,n} = 0.
$$ 
This implies that $x_{m,n}=x_m$ is independent of the index $n$ while $y_{m,n} = n\, f(m) + g(m)$, with $\Delta x_m \, \delta y_{m,n} = (x_{m+1}-x_m)\cdot f(m)\neq 0$.  Making the substitutions \eqref{2D difference operators}, the expressions \eqref{derivative computational variable expressions} are approximated by 
\begin{equation}\label{discrete derivatives}
\begin{aligned}
&u_y^{m,n} = \frac{\delta u_{m,n}}{\delta y_{m,n}},& &
u_x^{m,n} = \frac{\delta u_{m,n}\, \Delta u_{m,n} - \Delta y_{m,n}\, \delta u_{m,n}}{\Delta x_m\, \delta y_{m,n}},\\
&u_{yy}^{m,n} = \frac{\delta^2 u_{m,n}}{(\delta y_{m,n})^2},& &
u_{xy}^{m,n} = \frac{\delta y_{m,n}\, \Delta \delta u_{m,n} - \Delta\delta y_{m,n}\, \Delta u_{m,n} - \Delta y_{m,n}\, \delta^2 u_{m,n}}{\Delta x_m\, (\delta y_{m,n})^2},\\
&u_{yyy}^{m,n} = \delta^3 u_{m,n}{(\delta y_{m,n})^3},&\quad &
u_{xyy}^{m,n} = \frac{\delta y_{m,n}\, \Delta \delta^2 u_{m,n} - 2 \delta^2u_{m,n}\, \Delta\delta y_{m,n} - \Delta y_{m,n}\, \delta^3 u_{m,n}}{\Delta x_m (\delta y_{m,n})^3}.
\end{aligned}
\end{equation}
\end{example}

We are now interested in the induced action of a Lie pseudo-group on discrete points.

\begin{definition}
Given a Lie pseudo-group $\G$ acting on $M$, the \emph{pseudo-group product} action on the $k$-fold Cartesian product $M^{\times k}$ is
\begin{equation}\label{product action}
(g\cdot z_{N_1}, \ldots, g\cdot z_{N_k}),\qquad g\in \G,
\end{equation}
provided the points $z_{N_1},\ldots,z_{N_k} \in \text{dom } g$.
\end{definition}

\begin{remark} 
The nature of the product action \eqref{product action} depends on the type of the Lie pseudo-group $\G$.  If the Lie pseudo-group $\G$ is of infinite type, its $k$-fold product action is no longer a Lie pseudo-group as it is not possible to encapsulate into a system of differential equations the requirement that the same diffeomorphism should act on different points.   In this case, the product action only satisfies the defining properties of a pseudo-group.  On the other hand, the $k$-fold product action of a Lie pseudo-group of finite type, i.e. a local Lie group action, remains a Lie pseudo-group of finite type.  Another important distinction between pseudo-groups of finite and infinite types occurs when more copies of the manifold $M$ are appended to the Cartesian product $M^{\times k}$.  
For pseudo-groups of infinite type, when a new copy of $M$ is added, new pseudo-group parameters will occur in the product action while this is not the case for pseudo-groups of finite type. 
\end{remark}

As shown in Example \ref{example intro}, no joint invariant of the product action \eqref{product action example 2} can approximate the differential invariant \eqref{I}. This is not peculiar to this pseudo-group and another instance is given in Example \ref{second example}.  To address this lack of joint invariants, we propose to discretize the product pseudo-group action, replacing derivatives by finite difference approximations.   Before stating the general theory, our proposed idea is applied to the product pseudo-group \eqref{product action example 2}.

\begin{example}\label{discretization action example}
On the rectangular grid
$$
\delta x_{m,n} = x_{m,n+1}-x_{m,n}=0,\qquad y_n= k\, n + y_0,
$$
a suitable discretization of the product pseudo-group action \eqref{product action example 2} is obtained by approximating the first order derivative $f_x(x_{m})$ by the forward difference
\begin{equation}\label{fprime discretization}
f_x(x_{m}) \sim f^m_x = \frac{\Delta f_m}{\Delta x_m} = \frac{f_{m+1} - f_m}{x_{m+1} - x_m}\quad \text{where}\quad f_{m+j}=f(x_{m+j}),
\end{equation}
to give the discretized action
\begin{equation}\label{approximate pseudo-group action}
\G_d \colon \qquad X_m=f_m,\qquad Y_n=y_n,\qquad U_{m,n} = \frac{u_{m,n}}{f^m_x}.
\end{equation}
The subscript $d$ is added to $\G$ to indicate that the pseudo-group action has been discretized.  For \eqref{approximate pseudo-group action} to be a legitimate discretization it must satisfy the properties of an action.  These are readily seen to be satisfied except maybe for closure under composition. To this end, let
$$
\widetilde{X}_m=\widetilde{f}_m=\widetilde{f}(X_m),\qquad \widetilde{Y}_n = Y_n= y_n,\qquad \widetilde{U}_{m,n}= \frac{U_{m,n}}{\widetilde{f}^m_X},
$$
be a second discretized transformation. Then $\widetilde{X}_{m}=\widetilde{f} \comp f(x_m)$, and
\begin{align*}
\widetilde{U}_{m,n} &=U_{m,n}\cdot \frac{\Delta X_m}{\Delta[ \widetilde{f}(X_m)]}
=  u_{m,n} \cdot \frac{\Delta x_m}{\Delta [f(x_m)]} \cdot \frac{\Delta [f(x_m)]}{\Delta [\widetilde{f} \comp f(x_m)]}\\
&=  u_{m,n} \cdot \frac{\Delta x_m}{\Delta [\widetilde{f} \comp f(x_m)]}
= \frac{u_{m,n}}{(\widetilde{f} \comp f)^m_x},
\end{align*}
showing that \eqref{approximate pseudo-group action} is closed under composition.  

The approximation \eqref{fprime discretization} is not unique.  Any other discretization preserving the group action properties is acceptable.  For example, the approximation \eqref{fprime discretization} could be replaced by the backward difference
$$
f_x(x_m) \sim \frac{f_m - f_{m-1}}{x_m - x_{m-1}}.
$$
On the other hand, \eqref{approximate pseudo-group action} is not closed under composition if the centred approximation
$$
f_x(x_m) \sim \frac{1}{2} \bigg[ \frac{\Delta f_{m}}{\Delta x_m} + \frac{\Delta f_{m-1}}{\Delta x_{m-1}} \bigg]
$$
is considered.

At the infinitesimal level, the discretized action \eqref{approximate pseudo-group action} is generated by the vector field
\begin{equation}\label{va}
\vv_a = a_m\pp{}{x_m} - u_{m,n}\, a_x^m \pp{}{u_{m,n}},
\end{equation}
where
$$
a_m=a(x_m)\qquad \text{and}\qquad a_x^m = \frac{a_{m+1} - a_m}{x_{m+1} - x_m}.
$$
To compute the Lie algebra structure of \eqref{va}, one would use the standard prolongation, \cite{LW-2006},
$$
\pr\vv_a = \sum_{m}\> a_m \pp{}{x_m} - \sum_{m,n}\> u_{m,n}\, a_x^m\pp{}{u_{m,n}}
$$
and define the Lie bracket $[\vv_a, \vv_b]$ of two infinitesimal generators to be the vector field satisfying 
$$
\pr[\vv_a,\vv_b]:= \pr \vv_a \comp \pr \vv_b - \pr \vv_b \comp \pr \vv_a.
$$
For two infinitesimal generators of the form \eqref{va}, we obtain the expected commutation relation
$$
[\vv_a, \vv_b] = \vv_{ab^\prime - ba^\prime}.
$$
\end{example}

\begin{remark}
By introducing the approximation \eqref{fprime discretization}, the discretized action \eqref{approximate pseudo-group action}  is no longer local as the approximation \eqref{fprime discretization} introduces the extra independent variable $x_{m+1}$ into the action at $(x_m, y_n, u_{m,n})$.  This type of non-local discrete transformations is reminiscent of transformations obtained when considering discrete generalized symmetries, \cite{LW-2006}.  Similar pseudo-group discretization also recently appeared in a discrete version of Noether's Second Theorem, \cite{HM-2011}.
\end{remark}




Given an admissible discrete pseudo-group action, it is possible to implement the moving frame method in a fashion similar to the continuous setting.  In the continuous theory, the jets of functions occurring in the prolonged action play the role of the \mbox{pseudo-group} parameters.  In the discrete case, the functions evaluated at distinct points will play the role of the \mbox{pseudo-group} parameters.  

\begin{example}
At the point $(x_m, y_n, u_{m,n})$, the discrete pseudo-group action
$$
X_m = f_m=f(x_m),\qquad Y_n = y_n,\qquad U_{m,n} = u_{m,n} \cdot \frac{x_{m+1}-x_m}{f_{m+1}-f_m}
$$
involves the pseudo-group parameters $g_{m,n} = (f_m,f_{m+1})$.
\end{example}

In the following, the pseudo-group parameters occurring in the discrete pseudo-group action at $z_N$ is denoted $g_N$.  To approximate the $n^\text{th}$ order lifted bundle 
\eqref{lifted bundle} we introduce the $n^{\text{th}}$ order \emph{forward joint lifted bundle} $\mathcal{B}^{[n]}$  parameterized by $(\mathfrak{z}^{[n]}_N, g^{[n]}_N)$, where
$$
g^{[n]}_N = (\ldots\, g_{N+K}\, \ldots),\qquad 0 \leq \# K \leq n.
$$
The fibre of the $n^\text{th}$ order forward joint lifted bundle $\mathcal{B}^{[n]}$ at $\mathfrak{z}^{[n]}_N$ is denoted $\G_N^{[n]}$.  The discretized pseudo-group $\G_d$ acts on $\mathcal{B}^{[n]}$ by right multiplication
$$
R_{h_N}(\mathfrak{z}^{[n]}_N, g^{[n]}_N) = (h_N^{[n]}\cdot \mathfrak{z}^{[n]}_N, (g\cdot h^{-1})^{[n]}|_{h_N^{[n]}\cdot \mathfrak{z}^{[n]}_N}).
$$

\begin{definition}
Let $\G_d$ be a discretized Lie pseudo-group acting on the $n^\text{th}$ order joint lifted bundle $\mathcal{B}^{[n]}$.  An order $n$ (right) \emph{joint moving frame} is a $\G_d$-equivariant section of the order $n$ joint lifted bundle $\mathcal{B}^{[n]}$:
$$
\widehat{\rho}^{\,[n]}(\mathfrak{z}^{[n]}_N) =  (\mathfrak{z}^{[n]}_N,\,\rho^{[n]}(\mathfrak{z}^{[n]}_N)).
$$
\end{definition}

\noindent Right equivariance means that
$$
R_{g_N}\widehat{\rho}^{\, [n]}(\mathfrak{z}^{[n]}_N) = \widehat{\rho}^{\,[n]}(g_N^{[n]}\cdot \mathfrak{z}^{[n]}_N). 
$$
As in the continuous setting, a moving frame exists on (an open set of) the $n^\text{th}$ order joint bundle $\mathcal{J}^{[n]}$ if the action is free and regular.

\begin{definition}
A discretized Lie pseudo-group $\G_d$ acts freely at $\mathfrak{z}^{[n]}_N$ if the isotropy group
\begin{equation}\label{freeness condition}
\G^{[n]}_{N;\mathfrak{z}^{[n]}_N}=\{ g^{[n]}_N\,:\, g^{[n]}_N \cdot \mathfrak{z}^{[n]}_N=\mathfrak{z}^{[n]}_N \} = \{ \mathds{1}^{[n]}_N\},
\end{equation}
where $\mathds{1}^{[n]}_N$ is the discrete identity transformation at $\mathfrak{z}_N^{[n]}$.
\end{definition} 

\begin{example}
For the discretized pseudo-group \eqref{approximate pseudo-group action}, the isotropy condition $\g_{m,n}^{[0]} \cdot \mathfrak{z}_{m,n}^{[0]} = \mathfrak{z}_{m,n}^{[0]}$ is
$$
x_m = f_m,\qquad y_n = y_n,\qquad u_{m,n} = \frac{u_{m,n}}{f_x^m}
$$
which requires
$$
f_m=x_m,\qquad f_{m+1}=x_{m+1}.
$$
In general, the isotropy condition $\g_{m,n}^{[k]} \cdot \mathfrak{z}_{m,n}^{[k]} = \mathfrak{z}^{[k]}_{m,n}$ yields
$$
f_{m+\ell} = x_{m+\ell},\qquad \ell = 0,\ldots,k+1.
$$
\end{example}


Provided the discrete product pseudo-group action is free and regular, a joint moving frame is constructed through a normalization procedure similar to the continuous case.  Let $\mathcal{K}^{[n]}=\{ z_{i_1}=c_1,\ldots, z_{i_{r_n}}=c_{r_n} \} \subset \mathcal{J}^{[n]}$ be a coordinate cross-section, then the corresponding joint moving frame $\widehat{\rho}^{[n]}$ is obtained by solving the normalization equations
$$
Z_{i_1}(\mathfrak{z}^{[n]}_N,g^{[n]}_N)=c_1,\qquad \ldots \qquad Z_{i_{r_n}}(\mathfrak{z}^{[n]}_N,g^{[n]}_N)=c_{r_n}
$$
for the pseudo-group parameters $g^{[n]}_N=\rho^{[n]}(\mathfrak{z}^{[n]}_N)$.  For $k >n$, the cross-sections are required to be compatible, that is $\pi^k_n(\mathcal{K}^{[k]}) = \mathcal{K}^{[n]}$. The corresponding moving frames are then compatible $\widehat{\pi}^k_n\comp \widehat{\rho}^{\,[n]} = \widehat{\rho}^{\,[n]} \comp \pi^k_n$.  Here $\widehat{\pi}^k_n\colon \mathcal{B}^{[k]} \to \mathcal{B}^{[n]}$ is the standard projection obtained by truncation.  A discrete analogue of Theorem \ref{persistence of freeness} also holds.

\begin{theorem}\label{discrete persistence of freeness}
Let $\G_d$ be the discretization of a Lie pseudo-group and assume it acts freely at $\mathfrak{z}^{[n]}_N$ for any $N \in \mathbb{Z}^p$.  Then for $k>n$,  $\G_d$ acts freely at $\mathfrak{z}^{[k]}_N$.
\end{theorem}

\begin{remark}
Before proving Theorem \ref{discrete persistence of freeness} in general, it is instructive to consider a low dimensional example.  In two dimensions, assume the discretized action is free at $\mathfrak{z}_{m,n}^{[2]}$.  Our goal is to show that it remains free at $\mathfrak{z}_{m,n}^{[3]}$.  In Figure \ref{order 3 jet}, the multi-indices contained in $\mathfrak{z}^{[3]}_{m,n}$ are displayed.  Figures \ref{order 2 discrete jets 1}--\ref{order 2 discrete jets 3} show that sitting inside $\mathfrak{z}^{[3]}_{m,n}$ are the order 2 discrete jets
$$
\mathfrak{z}^{[2]}_{m,n},\qquad \mathfrak{z}^{[2]}_{m+1,n},\qquad \mathfrak{z}^{[2]}_{m,n+1}.
$$
Since these three order 2 jets cover $\mathfrak{z}_{m,n}^{[3]}$,
\begin{subequations}\label{isomorphism}
\begin{equation}
\mathfrak{z}_{m,n}^{[3]} \simeq (\mathfrak{z}_{m,n}^{[2]},\mathfrak{z}_{m+1,n}^{[2]},\mathfrak{z}_{m,n+1}^{[2]}).
\end{equation}
Similarly, at the pseudo-group level
\begin{equation}
g_{m,n}^{[3]} \simeq (g_{m,n}^{[2]}, g_{m+1,n}^{[2]}, g_{m,n+1}^{[2]}).
\end{equation}
\end{subequations}
Next, a pseudo-group transformation $g_{m,n}^{[3]} \in \G_{m,n}^{[3]}$ keeps $\mathfrak{z}_{m,n}^{[3]}$ fixed if and only if  it keeps $\mathfrak{z}_{m,n}^{[2]}$, $\mathfrak{z}_{m+1,n}^{[2]}$, $\mathfrak{z}_{m,n+1}^{[2]}$ fixed simultaneously.  That is
\begin{equation}\label{G3}
\G^{[3]}_{m,n;\mathfrak{z}^{[3]}_{m,n}} = \G^{[3]}_{m,n;\mathfrak{z}^{[2]}_{m,n}} \cap \G^{[3]}_{m,n;\mathfrak{z}^{[2]}_{m+1,n}} \cap\G^{[3]}_{m,n;\mathfrak{z}^{[2]}_{m+1,n}},
\end{equation}
where
\begin{align*}
\G^{[3]}_{m,n;\mathfrak{z}^{[2]}_{m,n}} &= \{ g_{m,n}^{[3]}\, : \, g_{m,n}^{[3]} \cdot \mathfrak{z}_{m,n}^{[2]} = g_{m,n}^{[2]} \cdot \mathfrak{z}_{m,n}^{[2]} = \mathfrak{z}_{m,n}^{[2]} \}, \\
\G^{[3]}_{m,n;\mathfrak{z}^{[2]}_{m+1,n}} &= \{ g_{m,n}^{[3]}\, : \, g_{m,n}^{[3]} \cdot \mathfrak{z}_{m+1,n}^{[2]} = g_{m+1,n}^{[2]} \cdot \mathfrak{z}_{m+1,n}^{[2]} = \mathfrak{z}_{m+1,n}^{[2]} \},\\
\G^{[3]}_{m,n;\mathfrak{z}^{[2]}_{m,n+1}} &= \{ g_{m,n}^{[3]}\, : \, g_{m,n}^{[3]} \cdot \mathfrak{z}_{m,n+1}^{[2]} = g_{m,n+1}^{[2]} \cdot \mathfrak{z}_{m,n+1}^{[2]} = \mathfrak{z}_{m,n+1}^{[2]} \}.
\end{align*}
By assumption
$$
\G^{[2]}_{m,n;\mathfrak{z}^{[2]}_{m,n}} = \{ \mathds{1}_{m,n}^{[2]} \},\qquad \G^{[2]}_{m+1,n;\mathfrak{z}^{[2]}_{m+1,n}} = \{ \mathds{1}_{m+1,n}^{[2]} \},\qquad 
\G^{[2]}_{m,n+1;\mathfrak{z}^{[2]}_{m,n+1}} = \{ \mathds{1}_{m,n+1}^{[2]} \},  
$$
and it follows that
\begin{equation}\label{isotropy groups}
\begin{gathered}
\G^{[3]}_{m,n;\mathfrak{z}^{[2]}_{m,n}} \simeq \{ (\mathds{1}^{[2]}_{m,n}, *, *) \}, \qquad
\G^{[3]}_{m,n;\mathfrak{z}^{[2]}_{m+1,n}} \simeq \{ (*, \mathds{1}^{[2]}_{m+1,n},*) \},\\
\G^{[3]}_{m,n;\mathfrak{z}^{[2]}_{m,n+1}} \simeq \{ (*, *, \mathds{1}^{[2]}_{m,n+1}) \},
\end{gathered}
\end{equation}
under the isomorphism \eqref{isomorphism}.  The exact expressions for $*$ will depend on the particular pseudo-group action. Combining \eqref{G3} and \eqref{isotropy groups}, we conclude that
$$
\G^{[3]}_{m,n;\mathfrak{z}^{[3]}_{m,n}} \simeq \{ (\mathds{1}_{m,n}^{[2]}, \mathds{1}_{m+1,n}^{[2]}, \mathds{1}_{m,n+1}^{[2]})\} \simeq \{ \mathds{1}_{m,n}^{[3]} \}.
$$

\begin{figure}[!ht]
\setlength{\unitlength}{0.2mm}
\begin{center}
\subfloat[$\mathfrak{z}_{m,n}^{[3]}$]{
\begin{picture}(100,175)
\put(0,30){\vector(1,0){150}}
\put(10,20){\vector(0,1){150}}
\put(10,30){\circle*{10}}
\put(10,70){\circle*{10}}
\put(10,110){\circle*{10}}
\put(10, 150){\circle*{10}}
\put(50,30){\circle*{10}}
\put(50,70){\circle*{10}}
\put(50,110){\circle*{10}}
\put(90,30){\circle*{10}}
\put(90,70){\circle*{10}}
\put(130,30){\circle*{10}}
{\thicklines
\color{red}\put(10,30){\line(0,1){120}}
\color{red}\put(10,30){\line(1,0){120}}
\color{red}\put(130,30){\line(-1,1){120}}}\label{order 3 jet}
\end{picture}
}
\hskip 1.5cm
\subfloat[$\mathfrak{z}_{m,n}^{[2]}$]{
\begin{picture}(100,175)
\put(0,30){\vector(1,0){150}}
\put(10,20){\vector(0,1){150}}
\put(10,30){\circle*{10}}
\put(10,70){\circle*{10}}
\put(10,110){\circle*{10}}
\put(10, 150){\circle*{10}}
\put(50,30){\circle*{10}}
\put(50,70){\circle*{10}}
\put(50,110){\circle*{10}}
\put(90,30){\circle*{10}}
\put(90,70){\circle*{10}}
\put(130,30){\circle*{10}}
{\thicklines
\color{red}\put(10,30){\line(0,1){80}}
\color{red}\put(10,30){\line(1,0){80}}
\color{red}\put(90,30){\line(-1,1){80}}}
\label{order 2 discrete jets 1}
\end{picture}}
\hskip 1.5cm
\subfloat[$\mathfrak{z}_{m+1,n}^{[2]}$]{
\begin{picture}(100,110)
\put(0,30){\vector(1,0){150}}
\put(10,20){\vector(0,1){150}}
\put(10,30){\circle*{10}}
\put(10,70){\circle*{10}}
\put(10,110){\circle*{10}}
\put(10, 150){\circle*{10}}
\put(50,30){\circle*{10}}
\put(50,70){\circle*{10}}
\put(50,110){\circle*{10}}
\put(90,30){\circle*{10}}
\put(90,70){\circle*{10}}
\put(130,30){\circle*{10}}
{\thicklines
\color{red}\put(50,30){\line(0,1){80}}
\color{red}\put(50,30){\line(1,0){80}}
\color{red}\put(130,30){\line(-1,1){80}}}
\end{picture}}
\hskip 1.5cm
\subfloat[$\mathfrak{z}_{m,n+1}^{[2]}$]{
\begin{picture}(100,110)
\put(0,30){\vector(1,0){150}}
\put(10,20){\vector(0,1){150}}
\put(10,30){\circle*{10}}
\put(10,70){\circle*{10}}
\put(10,110){\circle*{10}}
\put(10, 150){\circle*{10}}
\put(50,30){\circle*{10}}
\put(50,70){\circle*{10}}
\put(50,110){\circle*{10}}
\put(90,30){\circle*{10}}
\put(90,70){\circle*{10}}
\put(130,30){\circle*{10}}
{\thicklines
\color{red}\put(10,70){\line(0,1){80}}
\color{red}\put(10,70){\line(1,0){80}}
\color{red}\put(90,70){\line(-1,1){80}}}
\label{order 2 discrete jets 3}
\end{picture}}
\end{center}
\caption{Forward discrete jets of order 2 contained in $\mathfrak{z}^{[3]}_{m,n}$.}
\label{order 3 discrete jets}
\end{figure}
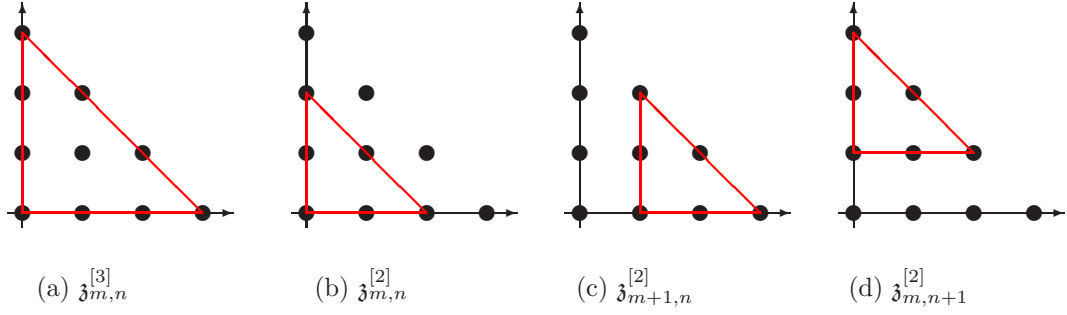
\end{remark}

\begin{proof}
To prove Theorem \ref{discrete persistence of freeness}, it suffices to consider the case when $k=n+1$ and proceed as in the 2-dimensional example above.  First, 
$$
\mathfrak{z}^{[n+1]}_N \simeq (\ldots\, \mathfrak{z}^{[n]}_{N+e_i}\, \ldots),\qquad i=1,\ldots,p.
$$
Next, since by assumption
$$
\G^{[n]}_{N+e_i;z^{[n]}_{N+e_i}} = \{ \mathds{1}^{[n]}_{N+e_i} \},\qquad i=1,\ldots,p,
$$
one has that 
\begin{align*}
\G^{[n+1]}_{N;\mathfrak{z}^{[n]}_{N+e_i}} &= \{ g_N^{[n+1]}\, :\, g_N^{[n+1]} \cdot \mathfrak{z}_{N+e_i}^{[n]} = g_{N+e_i}^{[n]} \cdot \mathfrak{z}_{N+e_i}^{[n]} =\mathfrak{z}_{N+e_i}^{[n]} \} \\
&\simeq \{ (\ldots,*, \ldots , \mathds{1}_{N+e_i}^{[n]},\ldots,*,\ldots)\},
\end{align*}
and
$$
\G^{[n+1]}_{N;\mathfrak{z}^{[n+1]}_{N}} = \bigcap_{i=1}^p\, \G^{[n+1]}_{N;\mathfrak{z}^{[n]}_{N+e_i}} 
= \{ \mathds{1}_{N}^{[n+1]} \}.
$$
\end{proof}


\begin{definition}
A function $I(\mathfrak{z}^{[n]}_N)\colon \mathcal{J}^{[n]} \to \R$ is a \emph{joint invariant} if
$$
I(g^{[n]}_N \cdot \mathfrak{z}^{[n]}_N) = I(\mathfrak{z}^{[n]}_N),\qquad g^{[n]}_N \in \G^{[n]}_N, 
$$ 
whenever the discrete product action is defined.
\end{definition}

\begin{definition}
Let $\widehat{\rho}^{\,[n]}(\mathfrak{z}^{[n]}_N)$ be an order $n$ joint moving frame.  The \emph{invariantization} of a function $F(\mathfrak{z}^{[n]}_N)$ is the joint invariant
\begin{equation}\label{discrete invariantization map}
\iota(F)(\mathfrak{z}^{[n]}_N) = F(\rho^{\,[n]}(\mathfrak{z}^{[n]}_N) \cdot \mathfrak{z}^{[n]}_N).
\end{equation}
\end{definition}

Of particular interest to us is the invariantization of the discrete derivatives $u^{\alpha;N}_{x^J}$:
\begin{equation}\label{joint invariants}
I^{\alpha;N}_{J} = \iota(u^{\alpha;N}_{x^J}),\qquad \alpha=1,\ldots,q,\qquad \# J \geq 0.
\end{equation}
We say that the cross-section $\mathcal{K}^{[n]}$ used to construct a joint moving frame $\widehat{\rho}^{[n]}$ is \emph{consistent} with the cross-section $\mathcal{K}\n$ used to construct a (continuous) moving frame $\widehat{\rho}\n$ if, in the continuous limit, $\mathcal{K}^{[n]}$ converges to $\mathcal{K}\n$.  For consistent cross-sections, since the discretized pseudo-group action $\G_d$ converges to the Lie pseudo-group $\G$ in the continuous limit, the discrete invariants \eqref{joint invariants} will converge to the differential invariants $I^\alpha_J = \iota(u^\alpha_{x^J})$: 
$$
\iota(u^{\alpha;N}_{x^J}) = I^{\alpha;N}_J \to I^\alpha_J = \iota(u^\alpha_{x^J}). 
$$ 

\begin{example}\label{joint moving frame 1st pseudo-group}
In this example, a joint moving frame for the discretized pseudo-group action \eqref{approximate pseudo-group action} is constructed.  First, a cross-section is given by
\begin{equation}\label{joint cross-section}
\mathcal{K}^{[\infty]} = \{ x_m = 0,\, u_{m+k,n}=1,\, k\in \mathbb{N}\}.
\end{equation}
Written differently, the cross-section is equivalent to 
\begin{equation}\label{equivalent cross-section}
x_m=0,\qquad u_{m,n}=1,\qquad \Delta^k(u_{m,n})=0,\qquad k=1,2, \ldots,
\end{equation}
which is an approximation of \eqref{cross-section} on a unit square mesh in the computational variables $(s,t)$.
Hence, \eqref{joint cross-section} is consistent with the cross-section \eqref{cross-section} used in the continuous setting.   Solving the normalization equations 
$$
0=X_m=f_m,\qquad 1=U_{m+k,n}=\frac{u_{m+k,n}}{f_x^{m+k}}=u_{m+k,n}\cdot\frac{\Delta x_{m+k}}{\Delta f_{m+k}},
$$
for the pseudo-group parameters $f_{m+k}$, $k\geq 0$, produces the (forward) joint moving frame
\begin{equation}\label{main pseudo-group joint moving frame}
f_m=0,\qquad f_{m+k} = \sum_{l=1}^{k} u_{m+l-1,n}\, \Delta x_{m+l-1},\qquad k=1, 2, 3, \ldots.
\end{equation}
Applying the invariantization map \eqref{discrete invariantization map} to the discrete variables $x_{m+k}$, $y_{n+l}$, $u_{m+k,n+l}$, we obtain the normalized joint invariants
\begin{equation}
\begin{aligned}
&\iota(x_{m+k})=(\rho^{[\infty]})^* f_{m+k}=
\begin{cases}
0 & \hskip 2.1cm k=0,\\
\displaystyle \sum_{l=1}^k u_{m+l-1,n}\, \Delta x_{m+l-1} &\hskip 2.1cm k=1,2, 3, \ldots,
\end{cases}\\
&\iota(y_{n+l})=y_{n+l},\qquad 
\iota(u_{m+k,n+l})=\frac{u_{m+k,n+l}}{(\rho^{[\infty]})^*f_x^{m+k}}=\frac{u_{m+k,n+l}}{u_{m+k,n}},\qquad k,l = 0, 1, 2, \ldots.
\end{aligned}
\end{equation}
Alternatively, invariantizing the forward differences in $x_m$, $y_n$, $u_{m,n}$ yields the joint invariants
\begin{gather*}
\iota(y_n)=y_n,\quad I^d_1= \iota(\Delta x_m)=u_{m,n}\, \Delta x_m,\quad \iota(\delta y_n)=k,\\
J_{0,1}^d=\iota(\delta u_{m,n}) = \frac{\delta u_{m,n}}{u_{m,n}},\quad I_2^d = \iota(\Delta^2 x_m) = \Delta u_{m,n}\, \Delta x_m + u_{m,n}\, \Delta^2 x_m,\quad \iota(\delta^2 y_n)=0,\\
J_{0,2}^d=\iota(\delta^2 u_{m,n}) = \frac{\delta^2 u_{m,n}}{u_{m,n}},\quad J_{1,1}^d=\iota(\Delta\delta u_{m,n}) =\frac{u_{m,n}\, \Delta\delta u_{m,n} - \delta u_{m,n}\, \Delta u_{m,n}}{u_{m+1,n}\, u_{m,n}}.
\end{gather*}
These invariants are finite difference approximations of the normalized differential invariants \eqref{normalized invariants pseudo-group1} on a unit square grid in the computational variables.
Another possibility is to invariantize the discrete derivatives
$$
u_{y}^{m,n} = \frac{\delta u_{m,n}}{\delta y_n},\qquad 
u_{xy}^{m,n}=\frac{\Delta \delta u_{m,n}}{\delta y_n \Delta x_m},\qquad
u_{yy}^{m,n} = \frac{\delta^2 u_{m,n}}{(\delta y_n)^2},
$$
to obtain the joint invariants
\begin{equation}\label{joint invariants 2}
\begin{gathered}
I_{0,1}^d=\iota(u^{m,n}_y)=\frac{u_y^{m,n}}{u_{m,n}},\qquad 
I_{0,2}^d=\iota(u_{yy}^{m,n})=\frac{u_{yy}^{m,n}}{u_{m,n}},\\
I_{1,1}^d=\iota(u_{xy}^{m,n})=\frac{u_{m+1,n+1}\, u_{m,n}-u_{m+1,n}\, u_{m,n+1}}{ u_{m,n}^2\,u_{m+1,n}\, \Delta x_m\,\delta y_n} = \frac{u_{m,n}\, u_{xy}^{m,n} - u_x^{m,n}\, u_y^{m,n}}{u_{m,n}^2\, u_{m+1,n}}.
\end{gathered}
\end{equation}
In the continuous limit, the invariants \eqref{joint invariants 2} converge to the normalized differential invariants obtained in \cite{OP-2008}.
\end{example}

%
%

Our main illustrative pseudo-group \eqref{main pseudo-group} was chosen for its simplicity.  This pseudo-group can be embedded into the larger pseudo-group 
\begin{equation}\label{larger pseudo-group}
X=f(x),\qquad Y=g(y),\qquad U = \frac{u}{f_x\, g_y},\qquad f,g \in \D(\mathbb{R}),
\end{equation}
which, as observed in the introduction, is the full symmetry group of the differential equation \eqref{main pde}.  In the following example, joint invariants of the discretized action \eqref{larger pseudo-group} are computed.  The results of these computations will be used in Sections \ref{numerical scheme section} and \ref{numerics section} to construct a fully invariant numerical scheme of equation \eqref{main pde} and perform numerical tests. 

\begin{example}\label{joint moving frame larger pseudo-group}
The construction of a joint moving frame for the pseudo-group \eqref{larger pseudo-group} is similar to the previous example.  Though, one important difference between the two examples is that it is no longer possible to work under the assumption that the step size $\delta y_n = k$ is constant as this is not an invariant constraint of the larger pseudo-group \eqref{larger pseudo-group}.  The most one can impose is that the mesh be rectangular
\begin{equation}\label{xy constraints}
\delta x_{m,n} = x_{m,n+1} - x_{m,n} = 0,\qquad \Delta y_{m,n} = y_{m+1,n} - y_{m,n}=0,
\end{equation}
so that $x_{m,n}=x_m$ and $y_{m,n}=y_n$.  At the discrete level, the pseudo-group action \eqref{larger pseudo-group} can be approximated by the forward discrete action
\begin{equation}\label{larger discretized pseudo-group}
X_m = f_m = f(x_m),\qquad Y_n = g_n = g(y_n),\qquad U_{m,n} = \frac{u_{m,n}}{f_x^m\, g_y^n},
\end{equation}
where
$$
f_x^m = \frac{\Delta f_m}{\Delta x_m} = \frac{f_{m+1}-f_m}{x_{m+1}-x_m},\qquad g_y^n = \frac{\delta g_n}{\delta y_n} = \frac{g_{n+1} - g_n}{y_{n+1}-y_n}.
$$
Summarizing the moving frame construction, a cross-section is given by
\begin{gather*}
\mathcal{K}^{[\infty]}\colon\qquad  
x_m=0,\qquad 
y_n=0,\qquad 
\Delta x_m\, \delta\Delta^2 u_{m,n}-\Delta^2 x_m\, \delta\Delta u_{m,n} =(\Delta x_m)^3\, \delta y_n,\\ 
u_{m+k,n}=u_{m,n+k}=1,\qquad  k \geq 0,
\end{gather*}
and the corresponding joint moving frame is
\begin{subequations}\label{larger pseudo-group moving frame}
\begin{equation}
\begin{aligned}
&f_m=0,& 
\qquad& 
f_{m+k} = \frac{\delta y_n}{g_{n+1}} \sum_{l=0}^{k-1} u_{m+l,n}\, \Delta x_{m+l},
\\
&g_n=0,&
\qquad &
g_{n+k} = \frac{g_{n+1}}{u_{m,n}\delta y_n} \sum_{l=0}^{k-1} u_{m,n+l}\, \delta y_{n+l},
\end{aligned}
\end{equation}
where $k \geq 1$, and
\begin{equation}
g_{n+1}=\cfrac{u_{m,n+1}\, u_{m,n}\, (\Delta x_m)^2\, (\delta y_n)^2}{
u_{m+1,n}\, \Delta x_{m+1} \Delta\bigg[ \cfrac{1}{u_{m,n}\, \Delta x_m} \Delta\bigg(\cfrac{u_{m,n+1}}{u_{m,n}}\bigg)\bigg]}.
\end{equation}

\end{subequations}
The applications of the invariantization map \eqref{discrete invariantization map} to the discrete variables $x_{m+k}$, $y_{n+l}$, $u_{m+k,n+l}$ yields the normalized joint invariants
\begin{align*}
&\iota(x_{m+k})=(\rho^{[\infty]})^* f_{m+k}=
\begin{cases}
0 & \hskip 1.55cm k=0,\\
\cfrac{\delta y_n}{g_{n+1}} \displaystyle\sum_{l=0}^{k-1} u_{m+l,n}\, \Delta x_{m+l} &\hskip 1.55cm k=1,2, 3, \ldots,
\end{cases}\\
&\iota(y_{n+k})=(\rho^{[\infty]})^* g_{n+k}=
\begin{cases}
0 & \hskip 1.5cm k=0,\\
\cfrac{g_{n+1}}{u_{m,n}\delta y_n} \displaystyle\sum_{l=0}^{k-1} u_{m,n+l}\, \delta y_{n+l} &\hskip 1.5cm k=1,2, 3, 4,\ldots,
\end{cases}\\ 
&\iota(u_{m+k,n+l})=\frac{u_{m+k,n+l}}{(\rho^{[\infty]})^*f_x^{m+k}\, (\rho^{[\infty]})^*g_y^{n+k}}=\frac{u_{m+k,n+l}\, u_{m,n}}{u_{m+k,n}\, u_{m,n+k}},\qquad k,l = 0, 1, 2, \ldots.
\end{align*}
For later use, the invariantization of  
$$
u_{xy}^{m,n} = \frac{\delta \Delta u_{m,n}}{\Delta x_m\, \delta y_n}
$$
gives the joint invariant
\begin{equation}\label{Id}
I_{1,1}^d=\iota(u_{xy}^{m,n})=\frac{u_{m+1,n+1}\, u_{m,n}-u_{m+1,n}\, u_{m,n+1}}{u_{m,n}\, u_{m+1,n}\, u_{m,n+1}\,\Delta x_m\,\delta y_n} = \frac{u_{m,n}\, u_{xy}^{m,n} - u_x^{m,n}\, u_y^{m,n}}{u_{m,n}\, u_{m+1,n}\, u_{m,n+1}}.
\end{equation}
\end{example}

\begin{example}\label{second example}
The Lie pseudo-group
\begin{equation}\label{2nd pseudo-group}
X = f(x),\qquad Y= e(x,y) = y\: f_x + g(x),\qquad U = u + \frac{e_x}{f_x} = u + \frac{y\, f_{xx} + g_x}{f_x},
\end{equation}
will serve as our last example.   This pseudo-group was used by Vessiot in his work on automorphic systems, \cite{V-1904}.  It is also one of the pseudo-groups used in \cite{OP-2008} to illustrate the method of equivariant moving frames.

By a similar argument to Example \ref{example intro}, on a generic mesh $(x_{m,n}, y_{m,n})$, the product pseudo-group action
\begin{equation}\label{product action 2nd pseudo-group}
\begin{gathered}
X_{m,n} = f(x_{m,n}),\qquad Y= e(x_{m,n},y_{m,n}) = y_{m,n}\: f_x(x_{m,n}) + g(x_{m,n}),\\
U_{m,n} = u_{m,n} + \frac{e_x(x_{m,n},y_{m,n})}{f_x(x_{m,n})} 
\end{gathered}
\end{equation}
has no joint invariant since $f(x_{m,n})$, $e(x_{m,n}, y_{m,n})$, and $e_x(x_{m,n},y_{m,n})/f_x(x_{m,n})$ are generically independent.  To reduce the number of pseudo-group parameters as much as possible, the invariant constraints
\begin{equation}\label{discrete x constraint}
\delta x_{m,n} = x_{m,n+1} - x_{m,n} = 0,\qquad \delta^2 y_{m,n} = y_{m,n+2}-2y_{m,n+1} + y_{m,n} = 0
\end{equation}
are imposed. Note that it is not possible to invariantly assume $\Delta y_{m,n}=y_{m+1,n} - y_{m,n} = 0$.  Hence, rectangular meshes are not invariant for this pseudo-group.  Provided $\delta y_{m,n}\neq 0$, which is an invariant constraint of \eqref{product action 2nd pseudo-group} when \eqref{discrete x constraint} is satisfied, the product pseudo-group action
$$
X_m = f(x_m),\quad Y= e(x_m,y_{m,n}) = y_{m,n}\: f_x(x_m) + g(x_m),\quad U_{m,n} = u_{m,n} + \frac{e_x(x_m,y_{m,n})}{f_x(x_m)} 
$$
admits the joint invariants
$$
\frac{y_{m,n+k}-y_{m,n}}{y_{m,n+1}-y_{m,n}},\qquad u_{m,n+k}-u_{m,n} - \bigg( \frac{y_{m,n+k}-y_{m,n}}{y_{m,n+1}-y_{m,n}} \bigg) (u_{m,n+1}-u_{m,n}),\qquad k \in \mathbb{Z}.
$$
By the same dilation argument as in Example \ref{example intro}, it is possible to conclude that these joint invariants cannot approximate all the differential invariants obtained in \cite{OP-2008}.  To construct further joint invariants the product action \eqref{product action 2nd pseudo-group} is discretized.  An admissible discretization is given by
\begin{equation}\label{pseudo-group discretization 2}
\begin{gathered}
X_m = f_m = f(x_m),\qquad Y_{m,n} = e_{m,n} = y_{m,n}\, f_{x}^m + g_m,\\
U_{m,n} =  u_{m,n} + \frac{e_{x}^{m,n}}{f_{x}^m} = u_{m,n} + \frac{\Delta e_{m,n}}{\Delta f_m}  - \frac{\Delta y_{m,n}}{\Delta x_m}=u_{m,n} + \frac{y_{m+1,n}}{\Delta f_{m}} \Delta\bigg( \frac{\Delta f_m}{\Delta x_m}\bigg) + \frac{\Delta g_m}{\Delta f_m},
\end{gathered}
\end{equation}
where
$$
g_m=g(x_m),\qquad f_{x}^m = \frac{\Delta f_m}{\Delta x_m}\qquad \text{and}\qquad e_{x}^{m,n} = \frac{\Delta e_{m,n}}{\Delta x_{m}} - \frac{\Delta y_{m,n}\, \delta e_{m,n}}{\Delta x_m\, \delta y_{m,n}}.
$$ 
To verify closure of \eqref{pseudo-group discretization 2} under composition, let
$$
\oX_m = \of_m,\qquad 
\oY_{m,n} = \overe_{m,n} = \frac{\Delta \of_m}{\Delta X_m}\, Y_{m,n} + \og_m,\qquad 
\oU_{m,n} = U_{m,n} + \frac{\Delta \overe_{mn}}{\Delta \of_m} - \frac{\Delta Y_{mn}}{\Delta X_m}.
$$
Thus, in the $x$ variable $\oX_m = \of_m = \of(X_m)=\overline{f}_m\comp f_m$, while in the $y$ variable
$$
\oY_{mn} = \overe_{m,n} = \frac{\Delta \of_m}{\Delta X_m}\, \bigg(  \frac{\Delta f_m}{\Delta x_m}\, y_{m,n} + g_m \bigg) + \og_m = \frac{\Delta \of_m}{\Delta x_m} \, y_{m,n} +  G_m,
$$
with
$$
G_m = \frac{\Delta \of_m}{\Delta f_m}\, g_m + \og_m.
$$
Finally, in the $u$ variable
\begin{align*}
\oU_{mn} &= u_{mn} + \frac{\Delta e_{mn}}{\Delta f_m} - \frac{\Delta y_{mn}}{\Delta x_m} + \frac{\Delta \overe_{mn}}{\Delta \of_m} - \frac{\Delta Y_{mn}}{\Delta X_m}\\
&=u_{mn} +  \frac{\Delta \overe_{mn}}{\Delta \of_m} - \frac{\Delta y_{mn}}{\Delta x_m},
\end{align*}
where the equality
$$
 \frac{\Delta e_{mn}}{\Delta f_m} = \frac{\Delta Y_{mn}}{\Delta X_m}
$$
was used.

To obtain a discrete approximation of the moving frame constructed in \cite{OP-2008}, we use the same cross-section replacing derivatives by their finite difference approximations:
\begin{gather*}
x_{m,n}=y_{m,n}=u_{m,n}=u_{x}^{m,n}=u_{y}^{m,n}=u_{xx}^{m,n}=u_{xy}^{m,n}=0,\\
 u_{yy}^{m,n}=1,\qquad u_{x^{k+3}}^{m,n}=u_{x^{k+2}y}^{m,n}=\cdots=0,\qquad k\geq 0.
\end{gather*}
The expressions for the discrete derivatives $u_y^{m,n}$, $u_{yy}^{m,n}$, $u_{yyy}^{m,n}$, and $u_{xyy}^{m,n}$ constraint to \eqref{discrete x constraint} appear in \eqref{discrete derivatives}.  Solving the normalization equations 
$$
X_{m,n}=Y_{m,n}=U_{m,n}=U_{X}^{m,n}=U_{Y}^{m,n}=U_{XY}^{m,n}=0,\qquad
U_{YY}^{m,n}=1,
$$
we obtain the pseudo-group normalizations
\begin{gather*}
f_m=0,\qquad e_{m,n}=0,\qquad \frac{\Delta f_m}{\Delta x_m} = \sqrt{u_{yy}^{m,n}},\qquad \frac{\Delta e_{m,n}}{\Delta x_m} = \sqrt{u_{yy}^{m,n}}\bigg(\frac{\Delta y_{m,n}}{\Delta x_m} - u_{m,n}\bigg),\\
\frac{\Delta f_{m+1}}{\Delta x_{m+1}} = \sqrt{u_{yy}^{m,n}} \bigg(1-\frac{\Delta x_m\, \delta u_{m,n}}{\delta u_{m+1,n}} \bigg),\qquad
\frac{\Delta e_{m+1}}{\Delta x_{m+1}} =\frac{\Delta f_{m+1}}{\Delta x_{m+1}} \bigg( \frac{\Delta y_{m+1,n}}{\Delta x_{m+1}} - u_{m+1,n}\bigg),\\
\frac{\Delta f_{m+2}}{\Delta x_{m+2}} = \frac{\Delta f_{m+1}}{\Delta x_{m+1}} \bigg(1 + \frac{\Delta x_{m+1}}{\delta y_{m+2,n}} [\delta y_{m,n}(\Delta y_{m,n}-\Delta x_m\, u_{m,n}) u_{yy}^{m,n} - u_{m+1,n+1}] \bigg).
\end{gather*}
The invariantization map \eqref{discrete invariantization map} provides the normalized joint invariants
\begin{gather*}
\iota(\Delta x_m)=\Delta x_m \sqrt{u_{yy}^{m,n}},\qquad
\iota(\Delta y_{m,n})=(\Delta y_{m,n} - u_{m,n}\, \Delta x_m)\sqrt{u_{yy}^{m,n}},\\
\iota(\delta y_{m,n})= \delta y_{m,n}\sqrt{u_{yy}^{m,n}},
\end{gather*}
and
\begin{equation}\label{joint invariants 2nd pseudo-group}
\begin{gathered}
I_{03}^d=\iota(u_{yyy}^{m,n})=\frac{u_{yyy}^{m,n}}{(u_{yy}^{m,n})^{3/2}},\qquad 
I_{1,2}^d=\iota(u_{xyy}^{m,n}) = \frac{u_{xyy}^{m,n} + u_{m,n}\, u_{yyy}^{m,n} + 2 u_{y}^{m,n}\, u_{yy}^{m,n}}{(u_{yy}^{m,n})^{3/2}}.
\end{gathered}
\end{equation}
In the continuous limit, the joint invariants \eqref{joint invariants 2nd pseudo-group} converge to the differential invariants 
$$
I_{0,3}^d \to I_{0,3}=\frac{u_{yyy}}{u_{yy}^{3/2}},\qquad I_{1,2}^d \to I_{1,2}=\frac{u_{xyy}+u\, u_{yyy}+2u_y\, u_{yy}}{u_{yy}^{3/2}},
$$
as obtained in \cite{OP-2008}.
\end{example}

\section{Differential and finite difference equations}\label{numerical scheme section}

This section recalls basic definitions pertaining to invariant differential equations and their invariant finite difference approximations, \cite{LW-2006,O-1993}.  To treat differential equations and finite difference equations on a similar footing, computational variables are introduced in the continuous setting.  Given a differential equation
\begin{equation}\label{differential equation}
\Delta(x,u\n)=0,
\end{equation}
the chain rule \eqref{uxJ} may be used to re-express \eqref{differential equation} in terms of $x^i=x^i(s)$, $u^\alpha=u^\alpha(s)$ and their computational derivatives $x^i_{s^A}$, $u^\alpha_{s^A}$:
\begin{subequations}\label{extended system}
\begin{equation}\label{computational differential equation}
\overline{\Delta}(s,x\n,u\n)=\Delta(x,u\n)=0,
\end{equation}
where $(x\n,u\n) = (\, \ldots\, x^i_{s^A}\, \ldots\, u^\alpha_{s^A}\, \ldots\,)$ on the left-hand side of \eqref{computational differential equation} and $u\n = (\, \ldots\, u^\alpha_{x^J}\, \ldots\,)$ on the right-hand side.  Equation \eqref{computational differential equation} can be supplemented by \emph{companion equations}, \cite{M-2010}, 
\begin{equation}\label{companion equations}
\widetilde{\Delta}(s,x\n,u\n)=0,
\end{equation}
\end{subequations}
which impose restrictions on the change of variables $s\mapsto x(s)$.  For the \emph{extended system} \eqref{extended system} to have the same solution space as the original equation \eqref{differential equation}, the companion equations \eqref{companion equations} cannot introduce differential constraints on the derivatives $u^\alpha_{s^A}$.  Also, they must respect the non-degeneracy condition \eqref{change of variables condition}.

\begin{definition}
A Lie pseudo-group $\G$ is said to be a \emph{symmetry (pseudo-)group} of a differential equation $\Delta(x,u\n)=0$ if for $g\in\G$,
$$
\Delta(g\cdot x, g\n\cdot u\n)=0\qquad \text{whenever}\qquad \Delta(x,u\n)=0.
$$
An extended system of differential equations $\{\overline{\Delta}(s,x\n,u\n)=0,\, \widetilde{\Delta}(s,x\n,u\n)=0\}$  is \emph{$\G$-compatible} with the $\G$-invariant differential equation $\Delta(x,u\n)=0$ if it is invariant under the pseudo-group $\G$:
$$
\begin{cases}\overline{\Delta}(s,g\n\cdot x\n,g\n\cdot u\n)=0,\\ \widetilde{\Delta}(s,g\n\cdot x\n,g\n\cdot u\n)=0,\end{cases} \qquad \text{whenever}\qquad 
\begin{cases}\overline{\Delta}(s,x\n,u\n)=0,\\ \widetilde{\Delta}(s,x\n,u\n)=0.\end{cases}
$$
\end{definition}


Using a perspective slightly different from the one introduced in \cite{BRW-2008,LW-2006}, a numerical scheme for the differential equation \eqref{differential equation}, or its extended counterpart \eqref{extended system}, is a set of finite difference equations
$$
E(\mathfrak{z}^{[n]}_N)=0,\qquad \widetilde{E}(\mathfrak{z}^{[n]}_N)=0,
$$
having the property that, in the continuous limit, these equations converge to the extended system \eqref{extended system}:
\begin{equation}\label{continuous limit}
E(\mathfrak{z}^{[n]}_N) \to \overline{\Delta}(s,x\n,u\n),\qquad \widetilde{E}(\mathfrak{z}^{[n]}_N)\to \widetilde{\Delta}(s,x\n,u\n).
\end{equation}

\begin{definition}
A discretized pseudo-group $\G_d$ is a \emph{symmetry group} of the numerical scheme $\{E(\mathfrak{z}^{[n]}_N)=0$, $\widetilde{E}(\mathfrak{z}^{[n]})=0\}$ if
$$
\begin{cases}E(g^{[n]}_N\cdot \mathfrak{z}^{[n]}_N)=0,\\ \widetilde{E}(g^{[n]}_N\cdot \mathfrak{z}^{[n]}_N)=0,\end{cases}
\qquad \text{whenever}\qquad 
\begin{cases}E(\mathfrak{z}^{[n]}_N)=0,\\ \widetilde{E}(\mathfrak{z}^{[n]}_N)=0.\end{cases}
$$
\end{definition}

Given a $\G$-invariant differential equation $\Delta = 0$, there are many different strategies for constructing an invariant numerical scheme, \cite{B-2013, KO-2004, LW-2006,O-2001,RV-2013}.  Assuming $\Delta$ is a differential invariant, one possibility is to obtain an invariant discretization $E$ of $\Delta$ using moving frames.  This can be done algorithmically by first discretizing $\Delta$ to obtain a finite difference approximation $F$.  Since this discretization is not necessarily invariant, see Example \ref{invariant scheme example} for an illustration of this fact, an invariant discretization of $\Delta$ is obtained by invariantizing $F$:
$$
\Delta \sim E = \iota(F).
$$
An invariant approximation of $\Delta =0$ is then given by $E=0$.  As for the mesh equations $\widetilde{E} = 0$, there is, unfortunately, no clear algorithm for determining these equations.  Nevertheless, there are obvious requirements that need to be satisfied.  First, these equations must include the invariant constraints occurring in the construction of a joint moving frame for the discretized pseudo-group action $\G_d$.  For example, in Example \ref{joint moving frame 1st pseudo-group}, the invariant constraint permitting the construction of a joint moving frame is given by 
\begin{equation}\label{x mesh constraint}
x_{m,n+1} - x_{m,n} = 0,
\end{equation}
and this equation would need to be part of the mesh equations of any invariant numerical scheme constructed from the joint moving frame \eqref{main pseudo-group joint moving frame}.  For some pseudo-group actions it might be possible to add further invariant mesh equations provided the non-degeneracy constraint \eqref{mesh non-degeneracy condition} is satisfied.  For example, in Example \ref{joint moving frame 1st pseudo-group}, since $Y_{m,n}=y_{m,n}$ is invariant, we assumed that $y_{n}= k\, n + y_0$ to simplify computations.  Under this assumption, equation \eqref{x mesh constraint} would be supplemented by the invariant mesh equations
$$
y_{m+1,n} - y_{m,n} = 0,\qquad y_{m,n+1} - y_{m,n} = k.
$$

\begin{example}\label{invariant scheme example}
A numerical scheme for the differential equation \eqref{main pde} invariant under the full (discretized) symmetry pseudo-group \eqref{larger discretized pseudo-group} is constructed.  Following the prescription above, the invariant \eqref{I} is naively discretized on a rectangular mesh
\begin{equation}\label{intermediate approximation}
I_{1,1} \sim F = \cfrac{u_{m+1,n+1}\,u_{m,n}-u_{m+1,n}\,u_{m,n+1}}{u_{m,n}^3\, \Delta x_m\, \delta y_n}.
\end{equation}
We note that this approximation is not invariant under \eqref{larger discretized pseudo-group}.  Using the results of Example \ref{joint moving frame larger pseudo-group}, an invariant approximation is obtained by invariantizing \eqref{intermediate approximation}:
$$
I_{1,1} \sim I_{1,1}^d = \iota(F) = \frac{u_{m+1,n+1}\, u_{m,n}-u_{m+1,n}\, u_{m,n+1}}{ u_{m,n}\, u_{m+1,n}\, u_{m,n+1}\,\Delta x_m\,\delta y_n}.
$$
The construction of the joint moving frame in Example \ref{joint moving frame larger pseudo-group} is based on the assumption that \eqref{xy constraints} holds. Hence, an invariant numerical scheme for \eqref{main pde} is given by
\begin{subequations}\label{invariant scheme}
\begin{equation}\label{inveq}
\frac{u_{m+1,n+1}\, u_{m,n}-u_{m+1,n}\, u_{m,n+1}}{ u_{m,n}\, u_{m+1,n}\, u_{m,n+1}\,\Delta x_m\,\delta y_n}=1,
\end{equation}
with mesh equations
\begin{equation}\label{mesh equations}
\delta x_{m,n}=0,\qquad
\Delta y_{m,n}=0.
\end{equation}
\end{subequations}
The scheme \eqref{invariant scheme} is an approximation of the differential equations
%
\begin{subequations}\label{st main pde}
\begin{equation}\label{main pde in computational variables}
I_{1,1}=\frac{u\, u_{st} - u_s\, u_t}{u^3\, x_s\, y_t} = 1
\end{equation}
and
\begin{equation}\label{main pde companion equations}
x_t=0,\qquad y_s=0,
\end{equation}
\end{subequations}
in the computational variables $x=x(s,t)$, $y=y(s,t)$, $u=u(s,t)$.  Equation \eqref{main pde in computational variables} is simply \eqref{main pde} expressed in computational variables while equations \eqref{main pde companion equations} are the invariant companion equations of the extended system \eqref{st main pde}.
\end{example}

\section{Numerical simulations}\label{numerics section}

In this section, the fully invariant numerical scheme \eqref{invariant scheme} is compared with the standard finite difference approximation
\begin{equation}\label{stdeq}
\begin{gathered}
\cfrac{u_{m+1,n+1}\,u_{m,n}-u_{m+1,n}\,u_{m,n+1}}{u_{m,n}^3\, \Delta x_m\, \delta y_n}=1,\\
\Delta x_{m,n}=h, \qquad \delta x_{m,n}=0,\qquad \Delta y_{m,n}=0,\qquad \delta y_{m,n}=k
\end{gathered}
\end{equation}
of equation \eqref{main pde}.  Since the mesh equations \eqref{mesh equations} do not specify the step sizes $\Delta x_{m,n}$ and $\delta y_{m,n}$, the equations 
$$
\Delta x_{m,n}=h,\qquad \delta y_{m,n}=k
$$
are supplemented to compare the two schemes on the same footing.  In other words, the numerical schemes \eqref{inveq} and \eqref{stdeq} are both defined over the same rectangular mesh.

\subsection{Methodology}

Equations \eqref{stdeq} and \eqref{inveq} both relate the values of the solution $u$ at the four corners of a rectangle on the mesh. Given, the value of $u$ at three corners, the equations provide the value of $u$ at the remaining vertex. These equations are suited for initial value problems (IVPs). For example, the value of $u$ in the $xy$-plane can be calculated if initial conditions on $u$ are specified on two perpendicular axes. Though, in practice, one has to limit itself to a finite rectangular domain and the specification of $u$ on two of its sides will completely determine the solution on the rectangle.  Figure \ref{4pointIVP} illustrates the situation on a $4\times 4$ rectangle. At each step the value of $u$ at the blue dot is a function of the solution at the green dots.  Filling the rectangle from left to right and then from bottom up, the whole rectangle is covered.

\begin{figure}[!h]
\setlength{\unitlength}{0.2mm}
\begin{center}
\subfloat[1st step.]{
\begin{picture}(110,145)
{\thicklines
\put(0,15){\line(1,0){120}}
\put(0,15){\line(0,1){120}}
}
\put(120,15){\line(0,1){120}}
\put(0,135){\line(1,0){120}}
\color{green}\put(0,15){\circle*{10}}
\put(0,55){\circle*{10}}
\put(40,15){\circle*{10}}
\color{blue}\put(40,55){\circle*{10}}
\end{picture}
}
\hskip 1.5cm
\subfloat[2nd step.]{
\begin{picture}(110,145)
{\thicklines
\put(0,15){\line(1,0){120}}
\put(0,15){\line(0,1){120}}
}
\put(120,15){\line(0,1){120}}
\put(0,135){\line(1,0){120}}
\put(0,15){\circle*{10}}
\put(0,55){\circle*{10}}
\color{green}\put(80,15){\circle*{10}}
\put(40,55){\circle*{10}}
\put(40,15){\circle*{10}}
\color{blue}\put(80,55){\circle*{10}}
\end{picture}
}
\hskip 1.5cm
\subfloat[9th step.]{
\begin{picture}(110,145)
{\thicklines
\put(0,15){\line(1,0){120}}
\put(0,15){\line(0,1){120}}
}
\put(120,15){\line(0,1){120}}
\put(0,135){\line(1,0){120}}
\put(0,15){\circle*{10}}
\put(40,15){\circle*{10}}
\put(80,15){\circle*{10}}
\put(80,55){\circle*{10}}
\put(120,15){\circle*{10}}
\put(120,55){\circle*{10}}
\put(0,55){\circle*{10}}
\put(0,95){\circle*{10}}
\put(40,55){\circle*{10}}
\put(40,95){\circle*{10}}
\put(0,135){\circle*{10}}
\put(40,135){\circle*{10}}
\color{green}
\put(80,95){\circle*{10}}
\put(80,135){\circle*{10}}
\put(120,95){\circle*{10}}
\color{blue}\put(120,135){\circle*{10}}
\end{picture}
}
\end{center}
\caption{Initial value problem on a rectangle.}
\label{4pointIVP}
\end{figure}
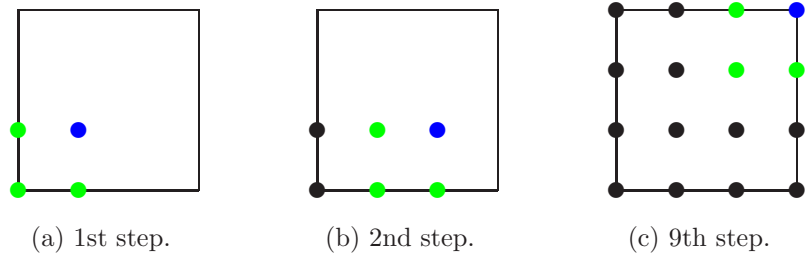

On the other hand, numerical schemes like \eqref{stdeq} and \eqref{inveq} are ill-defined for boundary value problems (BVP) on rectangular domains. Figure \ref{4pointBVP} illustrates the issue.  If, for example, one starts the iterative process in the bottom left corner of the domain of integration, then all points on the right and top boundaries highlighted in red in Figure \ref{4pointBVPend} are ill-defined since their values are simultaneously specified by the boundary conditions and the numerical scheme. 

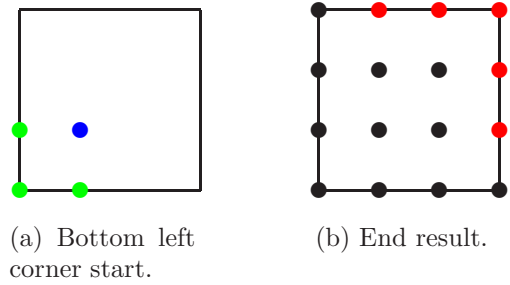
\begin{figure}[!h]
\setlength{\unitlength}{0.2mm}
\begin{center}
\subfloat[Bottom left corner start.]{
\begin{picture}(110,145)
{\thicklines
\put(0,15){\line(1,0){120}}
\put(0,15){\line(0,1){120}}
\put(120,15){\line(0,1){120}}
\put(0,135){\line(1,0){120}}
}\color{green}
\put(0,15){\circle*{10}}
\put(0,55){\circle*{10}}
\put(40,15){\circle*{10}}
\color{blue}\put(40,55){\circle*{10}}
\end{picture}
}
\hskip 1.5cm
\subfloat[End result.\label{4pointBVPend}]{
\begin{picture}(110,145)
{\thicklines
\put(0,15){\line(1,0){120}}
\put(0,15){\line(0,1){120}}
\put(120,15){\line(0,1){120}}
\put(0,135){\line(1,0){120}}
}
\put(0,15){\circle*{10}}
\put(40,15){\circle*{10}}
\put(80,15){\circle*{10}}
\put(80,55){\circle*{10}}
\put(0,55){\circle*{10}}
\put(0,95){\circle*{10}}
\put(40,55){\circle*{10}}
\put(40,95){\circle*{10}}
\put(80,95){\circle*{10}}
\put(0,135){\circle*{10}}
\put(120,15){\circle*{10}}
\color{red}\put(120,135){\circle*{10}}
\put(120,55){\circle*{10}}
\put(40,135){\circle*{10}}
\put(80,135){\circle*{10}}
\put(120,95){\circle*{10}}
\end{picture}
}
\end{center}
\caption{Ill-defined boundary value problem on a rectangle.}
\label{4pointBVP}
\end{figure}

Since in Section \ref{numerical subsection} we are interested in solving BVPs numerically, we now explain how to adapt the schemes \eqref{stdeq} and \eqref{inveq} to BVPs on rectangular domains. For this, we note that each point in the interior domain can be computed in four different ways using the numerical schemes. First, solving for $u_{m+1,n+1}$ in the invariant scheme \eqref{invariant scheme} we obtain
\begin{equation}\label{inv1}
u_{m+1,n+1}=u_{m+1,n}\,u_{m,n+1}\left(\cfrac{1}{u_{m,n}}+h\,k\right).
\end{equation}
Then, shifting \eqref{invariant scheme} from $(m,n)$ to $(m+1,n)$, the solution $u_{m+1,n+1}$ can also be expressed as
\begin{equation}\label{inv2}
u_{m+1,n+1}=\cfrac{u_{m+1,n}\,u_{m+2,n+1}}{u_{m+2,n}(1+h\,k\,u_{m+1,n})}.
\end{equation}
Similarly, shifting the invariant scheme \eqref{invariant scheme} from $(m,n)$ to $(m,n+1)$ and $(m+1,n+1)$ we obtain
\begin{gather}\label{inv3}
u_{m+1,n+1}=\cfrac{u_{m,n+1}\,u_{m+1,n+2}}{u_{m,n+2}(1+h\,k\,u_{m,n+1})},\qquad
u_{m+1,n+1}=\cfrac{u_{m+1,n}\,u_{m+2,n+1}}{u_{m+2,n}(1+h\,k\,u_{m+1,n})}.
\end{gather}
Defining $u_{m+1,n+1}$ to be the average of the four equations \eqref{inv1}, \eqref{inv2}, \eqref{inv3}  yields a finite difference equation expressing each interior point in the domain as a function of its eight surrounding points as illustrated in Figure \ref{BVP}.   The same procedure applies to the standard scheme \eqref{stdeq}.   The two new schemes are now well-adapted to BVPs on rectangular domains  since there is no conflict between the points computed using the numerical schemes and the boundary conditions.


\begin{figure}[h!]
\setlength{\unitlength}{0.2mm}
\begin{center}
\begin{picture}(200,215)
\put(0,15){$u_{m,n}$}
\put(80,15){$u_{m+1,n}$}
\put(180,15){$u_{m+2,n}$}
\put(180,130){$u_{m+2,n+1}$}
\put(0,130){$u_{m,n+1}$}
\put(0,210){$u_{m,n+2}$}
\put(80,210){$u_{m+1,n+2}$}
\put(80,130){$u_{m+1,n+1}$}
\put(180,210){$u_{m+2,n+2}$}
\color{green}
\put(40,35){\circle*{10}}
\put(120,35){\circle*{10}}
\put(200,35){\circle*{10}}
\put(40,115){\circle*{10}}
\put(40,195){\circle*{10}}
\put(200,195){\circle*{10}}
\put(120,195){\circle*{10}}
\put(200,115){\circle*{10}}
\color{blue}
\put(120,115){\circle*{10}}
\end{picture}
\caption{New scheme on nine points. The value of $u$ at the blue dot is determined by the neighbouring green points.}
\label{BVP}
\end{center}
\end{figure}
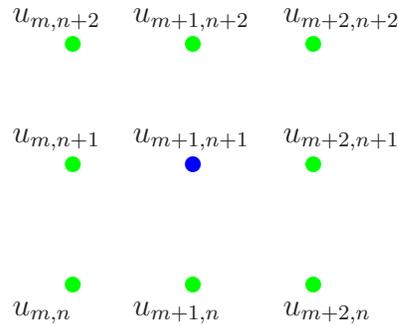
 
Solutions to BVPs are then obtained by applying the \emph{relaxation method}.  The first step in the implementation of the relaxation method consists of assigning values to the points inside the domain of integration.  In principle, arbitrary values can be assigned but it is always advantageous to assign well-educated initial values. In our case, we decided to use the average of the four solutions obtained by solving the IVPs starting in each corner of the rectangular domain.  Once this is done, new values are assigned to the interior points using the BVP adapted scheme.  Figure \ref{9point} illustrates the order in which one could assign these new interior values on a $5\times 5$ square.  Recomputing the interior values once using the most recent data is one iteration of the relaxation process. If the scheme is stable, by iterating the relaxation process the interior values will converge towards the scheme's solution.




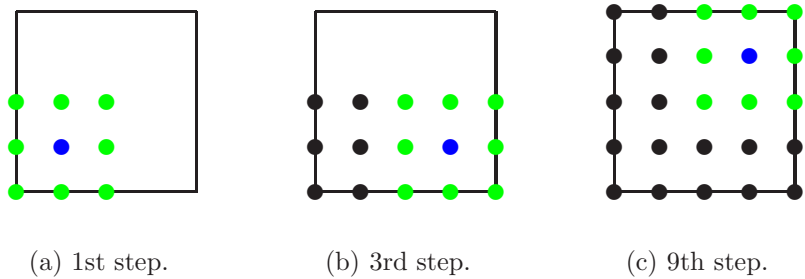
\begin{figure}[!h]
\setlength{\unitlength}{0.2mm}
\begin{center}
\subfloat[1st step.]{
\begin{picture}(110,175)
{\thicklines
\put(0,30){\line(1,0){120}}
\put(0,30){\line(0,1){120}}
\put(120,30){\line(0,1){120}}
\put(0,150){\line(1,0){120}}
}
\color{green}
\put(0,30){\circle*{10}}
\put(30,30){\circle*{10}}
\put(60,30){\circle*{10}}
\put(0,60){\circle*{10}}
\put(0,90){\circle*{10}}
\put(30,90){\circle*{10}}
\put(60,60){\circle*{10}}
\put(60,90){\circle*{10}}
\color{blue}\put(30,60){\circle*{10}}
\end{picture}
}
\hskip 1.5cm
\subfloat[3rd step.]{
\begin{picture}(110,175)
{\thicklines
\put(0,30){\line(1,0){120}}
\put(0,30){\line(0,1){120}}
\put(120,30){\line(0,1){120}}
\put(0,150){\line(1,0){120}}
}
\put(0,30){\circle*{10}}
\put(0,60){\circle*{10}}
\put(0,90){\circle*{10}}
\put(30,30){\circle*{10}}
\put(30,60){\circle*{10}}
\put(30,90){\circle*{10}}
\color{green}
\put(60,30){\circle*{10}}
\put(60,60){\circle*{10}}
\put(60,90){\circle*{10}}
\put(120,30){\circle*{10}}
\put(120,60){\circle*{10}}
\put(120,90){\circle*{10}}
\put(90,30){\circle*{10}}
\put(90,90){\circle*{10}}
\color{blue}\put(90,60){\circle*{10}}
\end{picture}
}
\hskip 1.5cm
\subfloat[9th step.]{
\begin{picture}(110,175)
{\thicklines
\put(0,30){\line(1,0){120}}
\put(0,30){\line(0,1){120}}
\put(120,30){\line(0,1){120}}
\put(0,150){\line(1,0){120}}
}
\put(0,30){\circle*{10}}
\put(0,60){\circle*{10}}
\put(0,90){\circle*{10}}
\put(30,30){\circle*{10}}
\put(30,60){\circle*{10}}
\put(30,90){\circle*{10}}
\put(60,30){\circle*{10}}
\put(60,60){\circle*{10}}
\put(120,30){\circle*{10}}
\put(120,60){\circle*{10}}
\put(90,30){\circle*{10}}
\put(90,60){\circle*{10}}
\put(0,120){\circle*{10}}
\put(0,150){\circle*{10}}
\put(30,120){\circle*{10}}
\put(30,150){\circle*{10}}
\color{green}
\put(60,90){\circle*{10}}
\put(90,90){\circle*{10}}
\put(120,90){\circle*{10}}
\put(60,150){\circle*{10}}
\put(90,150){\circle*{10}}
\put(120,150){\circle*{10}}
\put(120,120){\circle*{10}}
\put(60,120){\circle*{10}}
\color{blue}\put(90,120){\circle*{10}}
\end{picture}
}
\end{center}
\caption{Scheme on nine points covering a rectangular BVP.}
\label{9point}
\end{figure}


\subsection{Numerical results and analysis}\label{numerical subsection}

Three BVPs were tested using the exact solutions
\begin{gather}\label{exact solutions}
u=\cfrac{2}{(x+y)^2}, \qquad u=2\sec^2(x+y), 
\qquad u=\cfrac{2e^{x+y}}{(e^{x+y}-1)^2},
\end{gather}
obtained in \cite{P-2008}.
In each case, the boundary condition is given by the value of the exact solution on the edges of a rectangular domain. We note that the first and third solutions are not defined along the line $y+x=0$ and diverge to infinity on both sides of the singular line. The second solution also diverges along the lines $y+x=\pi/2+n\pi$, with $n\in\mathbb{Z}$. Since the quantitative results are similar for each solution, only the secant solution is presented below.

Table \ref{flat} lists the average error of the invariant and standard schemes \eqref{invariant scheme} and \eqref{stdeq} for different values of $h$ and $k$ for the secant solution on the unit square $[1,2]\times[1,2]$ after 100 iterations of the relaxation procedure.  For the cases considered, the invariant scheme is roughly three times more precise than the standard scheme.

\begin{table}[h!]
\begin{center}
\begin{tabular}{|*{5}{c|}}
\hline
Scheme&$h,k=0.1$&$h,k=0.05$&$h,k=0.01$&$h,k=0.005$\\
\hline
Standard&$2.19\times10^{-1}$&$1.07\times10^{-1}$&$3.23\times10^{-2}$&$1.66\times10^{-2}$\\
\hline
Invariant&$4.12\times10^{-2}$&$2.75\times10^{-2}$&$1.03\times10^{-2}$&$5.42\times10^{-3}$\\
\hline
\end{tabular}
\end{center}
\caption{Average errors on $[1,2]\times[1,2]$ for the secant solution after 100 relaxation iterations.}
\label{flat}
\end{table} 


As demonstrated in \cite{BRW-2008,KO-2004}, invariant schemes seem to shine near singularities. Here again the invariant scheme is more precise and stable near singularities. Table \ref{nearsing} shows the maximal error for both methods when the bottom left corner of the unit square of integration is brought closer to the exact solution singularity at $(\pi/4,\pi/4)\approx(0.785,0.785)$. The first row in the table gives the coordinates $(x_0,y_0)$ of the square of integration's bottom left corner. The step size in the independent variables is set to $h=k=0.01$, and the relaxation process was again run a hundred times.  As the square of integration gets closer to the singularity, Table \ref{nearsing} shows that the precision of the standard method gets worst much faster than the invariant scheme. Moreover, when $x_0=y_0=0.84$ or anywhere closer to the singularity $(\pi/4,\pi/4)$, the standard scheme becomes unstable while the invariant method integrates further into the singularity. As shown in Figure \ref{fig2}, while the source of the instability is in the bottom left corner, its manifestation appears first in the opposite corner for the standard method. Meanwhile, the invariant method, Figure \ref{fig3}, is faithful to the exact solution, Figure \ref{fig1}.

\begin{table}[h!]
\begin{center}
\begin{tabular}{|*{5}{c|}}
\hline
Scheme&$x_0=y_0=0.87$&$0.86$&$0.85$&$0.84$\\
\hline
Standard&$2.20$&$4.92$&$129.03$&$\text{unstable}$\\
\hline
Invariant&$3.12\times10^{-1}$&$4.46\times10^{-1}$&$6.78\times10^{-1}$&$1.15$\\
\hline
\end{tabular}
\end{center}
\caption{Maximal errors on a unit square near the singularity $(\pi/4,\pi/4)$.}
\label{nearsing}
\end{table}

\begin{figure}[!ht]
\psfrag{x}{$x$}
\psfrag{y}{$y$}
\psfrag{u(x,y)}{$u(x,y)$}
\begin{center}
\subfloat[Exact solution.\label{fig1}]{
\includegraphics[scale=0.5]{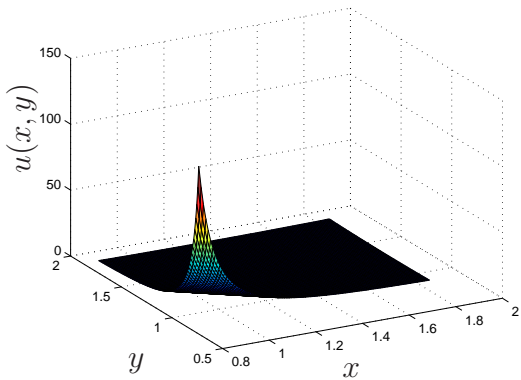}
}
\hskip 1cm
\subfloat[Standard scheme.\label{fig2}]{
\includegraphics[scale=0.5]{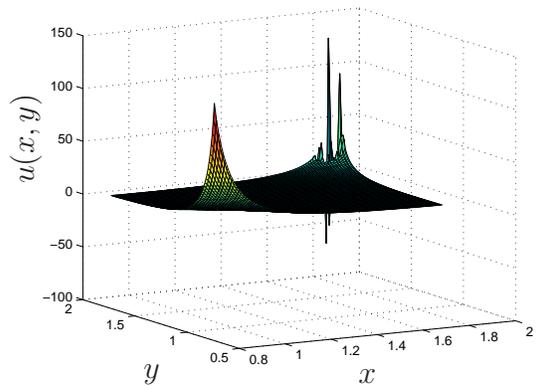}
}
\\
\subfloat[Invariant scheme.\label{fig3}]{
\includegraphics[scale=0.5]{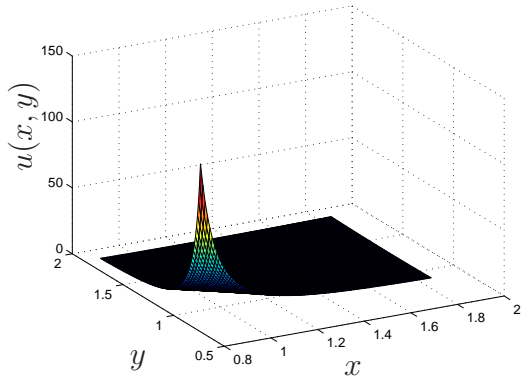}
}
\end{center}
\caption{Secant solution near the singularity $(\pi/4,\pi/4)$.}
\end{figure}

It is not difficult to understand why the invariant scheme produces better results when compared to the standard scheme.  The distinctive feature between the two schemes is the way the cubic term $u^3$ in \eqref{main pde} is approximated.  In the naive discretization \eqref{stdeq}, $u^3$ is approximated by the nonlinear term $u_{m,n}^3$. This cubic term in the standard scheme requires the use of a nonlinear equation solver like Newton's method at each iteration of the relaxation method which increases the computational cost and adds instability.  On the other hand, in the invariant scheme the cubic term $u^3$ is approximated by $u_{m,n}\,u_{m+1,n}\,u_{m,n+1}$.  By using the values of $u$ at three distinct points, the invariant method is more precise and stable, especially where the solution varies a lot. Moreover, \eqref{invariant scheme} can be solved for any of the $u$'s without the need to resort to a nonlinear solver. Thanks to this simplification, the computation time for the invariant method was approximately three times shorter than that of the standard method in all our numerical simulations.

As previously mentioned, similar results were also obtained for the rational and exponential solutions of \eqref{exact solutions}. 

\section{Conclusion}

To the best of our knowledge, this is the first work attempting to construct invariant numerical schemes of differential equations with infinite-dimensional symmetry groups.  As our examples show, the main issue with considering the product action of Lie pseudo-groups is the shortage of joint invariants to approximate differential invariants.  To circumvent this problem, we proposed to discretize the action by replacing derivatives with finite difference approximations.   To illustrate our constructions as clearly as possible, we chose simple Lie pseudo-group actions that have been well studied in the continuous setting.  The next natural step in this line of research would be to consider more substantial symmetry pseudo-groups and apply Lie theoretical tools to these invariant schemes to find explicit solutions.

The main emphasis of the paper was on the theoretical aspects that emerge when infinite-dimensional symmetry groups are discretized. Although the numerical simulations performed in Section \ref{numerics section} do not have the pretension to be the state of the art in numerical analysis, they indicate that invariant schemes can produce good numerical results.  It remains a challenge to bridge the gap between the most recent trends in numerical analysis and the latest developments in the theory of invariant finite difference equations.


\section*{Acknowledgement}

We thank Alexander Bihlo for stimulating discussions on the project, and Pavel Winternitz for his comments on the manuscript. The research of Rapha\"el Rebelo was supported in part by an FQRNT Doctoral Research Scholarship while the research of Francis Valiquette was supported in part by an AARMS Postdoctoral Fellowship.





\begin{thebibliography}{99}

\bibitem{B-2013}
Bihlo, A., Invariant meshless discretization schemes, {\it J. Phys. A: Math. Theor.} {\bf 46}: 6 (2013) 062001.

\bibitem{BRW-2008}
Bourlioux, A., Rebelo, R., and Winternitz, P., Symmetry preserving discretization of $SL(2,\mathbb{R})$ invariant equations, {\it J. Nonlin. Math. Phys.} {\bf 15}: 3 (2008) 362--372.

\bibitem{C-1953}
Cartan, \'E., Sur la structure des groupes infinis de transformations, in {\it Oeuvres Compl\`etes}, Part. II, Vol. 2, Gauthier--Villars, Paris, 1953, pp. 571--714.

\bibitem{CW-1987}
Champagne, C., and Winternitz, P., On the infinite-dimensional symmetry group of the Davey--Stewartson equations, {\it J. Math. Phys.} {\bf 29}: 1 (1988) 1--8.

\bibitem{DKLW-1985}
David, D., Kamran, N., Levi, D., and Winternitz, P.,  Subalgebras of loop algebras and symmetries of the Kadomtsev--Petviashvili equation, {\it Phys. Rev. Lett.} {\bf 55} (1985) 2111--2113.

\bibitem{D-2008}
Dawes, A.S.,  Invariant numerical methods, {\it Int. J. Numer. Meth. Fluids} {\bf 56} (2008) 1185--1191.


\bibitem{D-2011}
Dorodnitsyn, V.A., {\it Applications of Lie Groups to Difference Equations}, CRC Press, New York, 2011.

\bibitem{DKW-2003}
Dorodnitsyn, V., Kozlov, R., and Winternitz, P., Symmetry, Lagrangian formalism and integration of second order ordinary difference equations, {\it J. Nonlin. Math. Phys.} {\bf 10}: 2 (2003) 41--56.

\bibitem{HR-2011}
Huang, W., and Russell, R.D., {\it Adaptive Moving Mesh Methods}, Springer, New York, 2011.

\bibitem{IOV-2011}
Itskov, V., Olver, P.J., and Valiquette, F.,  Lie completion of pseudo-groups, {\it Transformation Groups} {\bf 16}:  1 (2011) 161--173.

\bibitem{HM-2011}
Hydon, P.E., and Mansfield, E.L., Extensions of Noether's second theorem:  from continuous to discrete systems, {\it Proc. R. Soc. A} {\bf 467} (2011) 3206--3221.

\bibitem{KO-2004}
Kim, P., and Olver, P.J., Geometric integration via multi-space, {\it Regular and Chaotic Dynamics} {\bf 9} (2004) 213--226.

\bibitem{K-1975}
Kumpera, A.,  Invariant diff\'erentiels d'un pseudogroupe de Lie, {\it J. Diff. Geom.} {\bf 10} (1975) 289--416.

\bibitem{K-1959}
Kuranishi, M., On the local theory of continuous infinite pseudo groups I, {\it Nagoya Math. J.} {\bf 15} (1959) 225--260.

\bibitem{L-1895}
Lie, S., Zur allgemeinen Theorie der partiellen Differentialgleichungen beliebeger Ordnung, {\it Leipz. Berich.} {\bf 47} (1895) 53--128; also {\it Gesammelte Abhandlungen}, Vol. 4, B.G. Teubner, Leipzig, 1929, 320--384.

\bibitem{LW-2006}
Levi, D., and Winternitz, P., Continuous symmetries of difference equations, {\it J. Phys. A: Math. Gen.} {\bf 39} (2006) R1--R63.


\bibitem{M-2010}
Mansfield, E.L., {\it A Practical Guide to the Invariant Calculus}, Cambridge University Press, Cambridge, 2010.


\bibitem{O-1993}
Olver, P.J., {\it Applications of Lie Groups to Differential Equations}, Springer--Verlag, New York, 1993.

\bibitem{O-1995}
Olver, P.J., {\it Equivalence, Invariants, and Symmetry}, Cambridge University Press, Cambridge, 1995.

\bibitem{O-2001}
Olver, P.J., Geometric foundations of numerical algorithms and symmetry, {\it Appl. Alg. Engin. Comp. Commun.} {\bf 11} (2001) 417--436.

\bibitem{OP-2005}
Olver, P.J., and Pohjanpelto, J.,  Maurer--Cartan forms and the structure of Lie pseudo-groups, {\it Selecta Math.} {\bf 11} (2005) 99--126.

\bibitem{OP-2008}
Olver, P.J., and Pohjanpelto, J., Moving frames for Lie pseudo-groups, {\it Canadian J. Math.} {\bf 60} (2008) 1336--1386.

\bibitem{OP-2009-1}
Olver, P.J., and Pohjanpelto, J., Differential invariant algebras of Lie pseudo-groups, {\it Adv. Math.} {\bf 222} (2009) 1746--1792.

\bibitem{OP-2009-2}
Olver, P.J., and Pohjanpelto, J., Persistence of freeness for Lie pseudogroup actions, Arkiv. Mat.  {\bf 50} (2012) 165--182.

\bibitem{P-2008}
Pohjanpelto, J., Reduction of exterior differential systems with infinite dimensional symmetry groups, {\it BIT Num. Math.} {\bf 48} (2008) 337--355.

\bibitem{RV-2013}
Rebelo, R., and Valiquette, F., Symmetry preserving numerical schemes for partial differential equations and their numerical tests, {\it J. Difference Eq. Appl.}, {\bf 19}: 5 (2013) 738--757.

\bibitem{RW-2004}
Rodr\'{\i}guez, M.A., and Winternitz, P., Lie symmetries and exact solutions of first-order difference schemes, {\it J. Phys. A: Math. Gen.} {\bf 37} (2004) 6129--6142.


\bibitem{SS-1965}
Singer, I.M., and Sternberg, S., The infinite groups of Lie and Cartan, part I, (The transitive groups), {\it J. Anal. Math.} {\bf 15} (1965) 1--114.


\bibitem{V-1904}
Vessiot, E., Sur l'int\'egration des syst\`emes diff\'erentiels qui admettent des groupes continus de transformations, {\it Acta. Math.} {\bf 28} (1904) 307--349.


\end{thebibliography}
\end{document}